\documentclass[11pt]{article}
\usepackage{latexsym,amssymb}
\usepackage{amsmath}
\usepackage{color}

\newtheorem{Theorem}{\bf Theorem}[section]
\newtheorem{Lemma}{\bf Lemma}[section]
\newtheorem{Proposition}{\bf Proposition}[section]
\newtheorem{Corollary}{\bf Corollary}[section]
\newtheorem{Remark}{\bf Remark}[section]
\newtheorem{Example}{\bf Example}[section]
\newtheorem{Definition}{\bf Definition}[section]

\newenvironment{theorem}{\begin{Theorem}$\!\!\!$}{\end{Theorem}}
\newenvironment{lemma}{\begin{Lemma}$\!\!\!$}{\end{Lemma}}

\newenvironment{corollary}{\begin{Corollary}$\!\!\!$}{\end{Corollary}}
\newenvironment{remark}{\begin{Remark}$\!\!\!$}{\end{Remark}}

\newenvironment{definition}{\begin{Definition}$\!\!\!$}{\end{Definition}}

\numberwithin{equation}{section}

\def\Xint#1{\mathchoice
{\XXint\displaystyle\textstyle{#1}}%
{\XXint\textstyle\scriptstyle{#1}}%
{\XXint\scriptstyle\scriptscriptstyle{#1}}%
{\XXint\scriptscriptstyle\scriptscriptstyle{#1}}%
\!\int}
\def\XXint#1#2#3{{\setbox0=\hbox{$#1{#2#3}{\int}$}
\vcenter{\hbox{$#2#3$}}\kern-.5\wd0}}

\def\dashint{\Xint-}

\begin{document}

\title{Existence of solutions for a fractional semilinear\\ parabolic equation 
with singular initial data}
\author{Kotaro Hisa and Kazuhiro Ishige}
\date{}
\maketitle
\begin{abstract}
In this paper we obtain necessary conditions and sufficient conditions 
on the initial data for the solvability of the Cauchy problem 
$$
\partial_t u+(-\Delta)^{\frac{\theta}{2}}u=u^p,\quad x\in{\bf R}^N,\,\,t>0,
\qquad
u(0)=\mu\ge 0\quad\mbox{in}\quad{\bf R}^N,
$$
where $N\ge 1$, $0<\theta\le 2$, $p>1$ and $\mu$ is a Radon measure 
or a measurable function in ${\bf R}^N$. 
Our conditions lead optimal estimates of the life span of the solution with 
$\mu$ behaving like $\lambda |x|^{-A}$ ($A>0$) at the space infinity, as $\lambda \to +0$. 
\end{abstract}
\section{Introduction}
This paper is concerned with the fractional semilinear parabolic equation 
\begin{equation}
\label{eq:1.1}
\partial_t u+(-\Delta)^{\frac{\theta}{2}}u=u^p,\qquad x\in{\bf R}^N,\,\,t>0,\\
\end{equation}
where $\partial_t:=\partial/\partial t$, $N\ge 1$, $0<\theta\le 2$ and $p>1$.
Here $(-\Delta)^{\theta/2}$ denotes the fractional power of the Laplace operator $-\Delta$ in ${\bf R}^N$. 
In this paper we show that 
every nonnegative solution of \eqref{eq:1.1} has a unique Radon measure in ${\bf R}^N$ 
as the initial trace and study qualitative properties of the initial trace. 
Furthermore, we give sufficient conditions for the existence of the solution of 
Cauchy problem \eqref{eq:1.1} and obtain optimal estimates of the life span of the solution with small initial data. 
\vspace{3pt}

Let us consider the case $\theta=2$, that is 
the semilinear parabolic equation 
\begin{equation}
\label{eq:1.2}
\partial_t u-\Delta u=u^p,\quad x\in{\bf R}^N,\,\,t>0,
\qquad
u(0)=\mu\ge 0\quad\mbox{in}\quad{\bf R}^N,
\end{equation}
where $N\ge 1$, $p>1$ and $\mu$ is a Radon measure or a measurable function in ${\bf R}^N$.  
The solvability of Cauchy problem 
\eqref{eq:1.2} has been studied extensively by many mathematicians 
since the pioneering work due to Fujita~\cite{F}
(see, for example, \cite{QS}, which is a book including a good list of references for problem~\eqref{eq:1.2}). 
Among others, in 1985, Baras and Pierre \cite{BP} 
proved the following by the use of the capacity of potentials of Meyers~\cite{M}. 
\begin{itemize}
  \item[(a)] 
  Let $u$ be a nonnegative local-in-time solution of \eqref{eq:1.2},  
  where $\mu$ is a Radon measure in ${\bf R}^N$. 
  Then $\mu$ must satisfy the following:
  \begin{itemize}
  \item
  If $\displaystyle{1<p<p_*}$, then 
  $\displaystyle{\sup_{x\in{\bf R}^N}\mu(B(x,1))<\infty}$;\\
  \item
  If $\displaystyle{p=p_*}$, then 
  $\displaystyle{\sup_{x\in{\bf R}^N}\mu(B(x,\sigma))\le \gamma\left|\log\sigma\right|^{-\frac{N}{2}}}$
  for all small enough $\sigma>0$;\\
  \item
  If $\displaystyle{p>p_*}$, then 
  $\displaystyle{\sup_{x\in{\bf R}^N}\mu(B(x,\sigma))\le \gamma\sigma^{N-\frac{2}{p-1}}}$
  for all small enough $\sigma>0$. 
  \end{itemize}
  Here $p_*:=1+2/N$ and $\gamma$ is a constant depending only on $N$ and $p$. 
\end{itemize}
Then we can find a positive constant $c_1$ with the following property:
\begin{itemize}
  \item[(b)] 
  Problem~\eqref{eq:1.2} possesses no local-in-time solutions 
  if $\mu$ is a nonnegative measurable function in ${\bf R}^N$ satisfying 
  $$
  \begin{array}{ll}
  \mu(x)\ge c_1|x|^{-N}\displaystyle{\biggr[\log\left(e+\frac{1}{|x|}\right)\biggr]^{-\frac{N}{2}-1}}
  \quad & \mbox{for}\quad \displaystyle{p=p_*},\vspace{7pt}\\
  \mu(x)\ge c_1|x|^{-\frac{2}{p-1}}\quad & \mbox{for}\quad \displaystyle{p>p_*},
  \end{array}
  $$
  in a neighborhood of the origin. 
\end{itemize}
For related results, see e.g., \cite{AD, BK}. 
On the other hand, Takahashi~\cite{Takahashi} recently proved that, 
in the case $p\ge p_*$,  for any $\gamma>0$, 
Cauchy problem~\eqref{eq:1.2} possesses no local-in-time nonnegative solutions 
with some Radon measure $\mu$ satisfying 
$$
\sup_{x\in{\bf R}^N}\mu(B(x,\sigma))\le \gamma\sigma^{N-\frac{2}{p-1}}\biggr[\log\biggr(e+\frac{1}{\sigma}\biggr)\biggr]^{-\frac{1}{p-1}}
$$
for all $\sigma>0$. See \cite[Theorem~1, Proposition~1]{Takahashi}.  
See also Remark~\ref{Remark:1.2}. 

The local solvability of Cauchy problem \eqref{eq:1.2} 
has been studied in many papers 
(see e.g., \cite{AD, BK, BC, FI, IKK, IKS, KY, RS, SL1, SL2, Takahashi, W1, W2} and references therein). 
It is known that 
there exists a constant $c_2>0$ such that 
Cauchy problem~\eqref{eq:1.2} possesses a solution in ${\bf R}^N\times[0,\rho^2]$, 
where $\rho>0$, if $p>p_*$ and 
\begin{equation}
\label{eq:1.3}
\sup_{x\in{\bf R}^N}\|\mu\|_{L^{r,\infty}(B(x,\rho))}\le c_2
\quad\mbox{with}\quad r=\frac{N(p-1)}{2} 
 \end{equation}
(see \cite{IKS}). See also \cite{IKK, KY, RS}.
This implies that, if $p>p_*$ and
$$
0\le\mu(x)\le c|x|^{-\frac{2}{p-1}}\quad\mbox{in}\quad{\bf R}^N\quad\mbox{with small enough $c>0$},
$$ 
then \eqref{eq:1.3} holds for any $\rho>0$ and 
problem~\eqref{eq:1.2} possesses a global-in-time solution. 
On the other hand, 
in the case $p=p_*$, 
as far as we know, 
there are no results on the local solvability 
of Cauchy problem~\eqref{eq:1.2} under such an assumption as 
\begin{equation}
\label{eq:1.4}
0\le\mu(x)\le c_3|x|^{-N}\displaystyle{\biggr[\log\left(e+\frac{1}{|x|}\right)\biggr]^{-\frac{N}{2}-1}}
\mbox{ in ${\bf R}^N$ with enough small $c_3>0$.}
\end{equation}

Some of the results on the solvability of Cauchy problem~\eqref{eq:1.2} 
are available to fractional semilinear parabolic equations, 
however there are no results on necessary conditions such as assertion~(a). 
\vspace{5pt}

In this paper we show the existence and the uniqueness of the initial trace of the solution of \eqref{eq:1.1} 
and obtain a refinement of assertion~(a). 
Furthermore, we give sufficient conditions on the existence of the solution of 
\begin{equation}
\label{eq:1.5}
\partial_t u+(-\Delta)^{\frac{\theta}{2}}u=u^p,\quad x\in{\bf R}^N,\,\,t>0,
\qquad
u(0)=\mu\ge 0\quad\mbox{in}\quad{\bf R}^N,
\end{equation}
where $N\ge 1$, $0<\theta\le 2$, $p>1$ and $\mu$ is a Radon measure or a measurable function in ${\bf R}^N$.  
Even in the case $\theta=2$, our sufficient conditions are new and they ensure that 
Cauchy problem~\eqref{eq:1.2} with~\eqref{eq:1.4} possesses a local-in-time solution. 
In addition, as an application of our conditions,  
we obtain optimal estimates of the life span of the solution of \eqref{eq:1.5} with 
$\mu=\lambda\phi$ as $\lambda\to +0$ by use of the behavior of $\phi$ at the space infinity. 
\vspace{5pt}

We introduce some notation and formulate the definition of the solutions of \eqref{eq:1.1}. 
For any $x\in{\bf R}^N$ and $r>0$, let $B(x,r):=\{y\in{\bf R}^N\,:\,|x-y|<r\}$ 
and $|B(x,r)|$ the volume of $B(x,r)$. 
Furthermore, for any $L^1_{{\rm loc}}({\bf R}^N)$ function $f$, we set 
$$
\dashint_{B(x,r)}\,f(y)\,dy:=\frac{1}{|B(x,r)|}\int_{B(x,r)} f(y)\,dy.
$$
Let $G=G(x,t)$ be the fundamental solution of 
\begin{equation}
\label{eq:1.6}
\partial_t u+(-\Delta)^{\frac{\theta}{2}}u=0\quad\mbox{in}\quad{\bf R}^N\times(0,\infty),
\end{equation}
where $0<\theta\le 2$. 
\begin{definition}
\label{Definition:1.1}
Let $u$ be a nonnegative measurable function in ${\bf R}^N\times(0,T)$, 
where $0<T\le\infty$. 
\vspace{3pt}
\newline
{\rm (i)} We say that $u$ is a solution of \eqref{eq:1.1} in ${\bf R}^N\times(0,T)$ 
if $u$ satisfies 
$$
\infty>u(x,t)=\int_{{\bf R}^N}G(x-y,t-\tau)u(y,\tau)\,dy+\int_\tau^t\int_{{\bf R}^N}G(x-y,t-s)u(y,s)^p\,dy\,ds
$$
for almost all $x\in{\bf R}^N$ and $0<\tau<t<T$.
\vspace{3pt}
\newline
{\rm (ii)} Let $\mu$ be a Radon measure in ${\bf R}^N$. 
We say that $u$ is a solution of \eqref{eq:1.5} 
in ${\bf R}^N\times[0,T)$ if $u$ satisfies 
\begin{equation}
\label{eq:1.7}
\infty>u(x,t)=\int_{{\bf R}^N}G(x-y,t)\,d\mu(y)+\int_0^t\int_{{\bf R}^N}G(x-y,t-s)u(y,s)^p\,dy\,ds
\end{equation}
for almost all $x\in{\bf R}^N$ and $0<t<T$. 
If $u$ satisfies \eqref{eq:1.7} with $=$ replaced by $\ge$, 
then $u$ is said to be a supersolution of 
\eqref{eq:1.5} in ${\bf R}^N\times[0,T)$.  
\vspace{3pt}
\newline
{\rm (iii)} 
Let $u$ be a solution of \eqref{eq:1.5} in ${\bf R}^N\times[0,T)$. 
We say that $u$ is a minimal solution of \eqref{eq:1.5} in ${\bf R}^N\times[0,T)$ if
$$
u(x,t)\le v(x,t)
\quad
\mbox{for almost all $x\in{\bf R}^N$ and $0<t<T$}
$$ 
for any solution $v$ of \eqref{eq:1.5} in ${\bf R}^N\times[0,T)$. 
\end{definition}
\vspace{5pt}

Now we are ready to state the main results of this paper. 
In the first theorem we show the existence and the uniqueness of the initial trace 
of the solution of \eqref{eq:1.1} and obtain a refinement of assertion~(a). 
See also Lemma~\ref{Lemma:2.4}. 
\begin{theorem}
\label{Theorem:1.1}
Let $N\ge 1$, $0<\theta\le 2$ and $p>1$. 
Let $u$ be a solution of \eqref{eq:1.1} in ${\bf R}^N\times(0,T)$, 
where $0<T<\infty$. 
Then there exists a unique Radon measure $\mu $ such that 
\begin{equation}
\label{eq:1.8}
\underset{t\to+0}{\mbox{{\rm ess lim}}}
\int_{{\bf R}^N}u(y,t)\phi(y)\,dy=\int_{{\bf R}^N}\phi(y)\,d\mu(y)
\end{equation}
for all $\phi\in C_0({\bf R}^N)$. 
Furthermore, 
there exists $\gamma_1>0$ depending only on $N$, $\theta$ and $p$ 
such that 
\begin{itemize}
  \item[{\rm (1)}] 
  $\displaystyle{\sup_{x\in{\bf R}^N}\mu(B(x,T^\frac{1}{\theta}))\le \gamma_1\,T^{\frac{N}{\theta}-\frac{1}{p-1}}}$ 
  if $1<p<p_\theta$;
  \item[{\rm (2)}] 
  $\displaystyle{\sup_{x\in{\bf R}^N}\mu(B(x,\sigma))\le \gamma_1\,
  \biggr[\log\biggr(e+\frac{T^{\frac{1}{\theta}}}{\sigma}\biggr)\biggr]^{-\frac{N}{\theta}}}$ 
  for all $0<\sigma<T^{\frac{1}{\theta}}$ if $p=p_\theta$; 
  \item[{\rm (3)}] 
  $\displaystyle{\sup_{x\in{\bf R}^N}\mu(B(x,\sigma))\le \gamma_1\,
  \sigma^{N-\frac{\theta}{p-1}}}$ 
  for all $0<\sigma<T^{\frac{1}{\theta}}$ if $p>p_\theta$.
\end{itemize} 
Here $p_\theta:=1+\theta/N$. 
\end{theorem}
\begin{remark}
\label{Remark:1.1}
{\rm (i)} Sugitani~{\rm\cite{Sugitani}} showed that, if $1<p\le p_\theta$ and $\mu\not\equiv 0$ in ${\bf R}^N$, 
then problem~\eqref{eq:1.5} possesses no nonnegative global-in-time solutions. 
\vspace{3pt}
\newline
{\rm (ii)} 
Let $u$ be a solution of \eqref{eq:1.1} in ${\bf R}^N\times[0,\infty)$
and $1<p\le p_\theta$.  
It follows from assertions~{\rm (1)} and {\rm (2)} that 
the initial trace of $u$ must be identically zero in ${\bf R}^N$. 
Then Theorem~{\rm\ref{Theorem:1.1}} 
leads the same conclusion as in Remark~{\rm\ref{Remark:1.1}}~{\rm (i)}. 
\end{remark}
As a corollary of Theorem~\ref{Theorem:1.1}, we have
\begin{corollary}
\label{Corollary:1.1}
Let $N\ge 1$, $0<\theta\le 2$ and $p>1$. 
Let $u$ be a solution of \eqref{eq:1.1} in ${\bf R}^N\times(0,T)$, where $0<T<\infty$.
Then there exists $\gamma>0$ depending only on $N$, $\theta$ and $p$ such that 
$$
\sup_{x\in{\bf R}^N}\,\dashint_{B(x,(T-t)^{\frac{1}{\theta}})} u(y,t)\,dy\le\gamma(T-t)^{-\frac{1}{p-1}}
$$
for almost all $0<t<T$. 
\end{corollary}
Our argument in the proof of Theorem~\ref{Theorem:1.1} is completely different 
from those in \cite{AD, BK, BP}. 
Let $u$ be a solution of \eqref{eq:1.1} in ${\bf R}^N\times(0,T)$, where $0<T<\infty$. 
We first prove the existence and the uniqueness of the initial trace of the solution~$u$. 
Next, in the case $p\not=p_\theta$ 
we apply the iteration argument in \cite[Theorem~5]{W1} 
to obtain an $L^\infty({\bf R}^N)$ estimate of the solution~$u$ (see Lemma~\ref{Lemma:3.1}). 
This yields a uniform estimate of $\|u(\tau)\|_{L^1(B(z,\rho))}$ with respect to $z\in{\bf R}^N$ and $\tau\in(0,T/2)$ 
for all small enough $\rho>0$ (see \eqref{eq:3.9}), 
and we complete the proof of Theorem~\ref{Theorem:1.1}. 
In the case $p=p_\theta$ 
we follow the argument in \cite{F, GK, Sugitani} 
and obtain an inequality related to  
$$
\int_{{\bf R}^N}u(x,t)G(x,t)\,dx
$$
(see \eqref{eq:3.22}). 
Then, applying the iteration argument in \cite[Section~2]{LN}, 
we prove Theorem~\ref{Theorem:1.1}. 
Furthermore, by Theorem~\ref{Theorem:1.1} we obtain
\begin{theorem}
\label{Theorem:1.2}
Assume the same conditions as in Theorem~{\rm\ref{Theorem:1.1}}.
Let $\mu$ be a Radon measure satisfying \eqref{eq:1.8}. Then $u$ is a solution of \eqref{eq:1.5} 
in ${\bf R}^N\times[0,T)$. 
\end{theorem}
We give sufficient conditions for the solvability of problem~\eqref{eq:1.5}. 
We modify the arguments in \cite{IKS, RS} and prove the following two theorems. 
\begin{theorem}
\label{Theorem:1.3}
Let $N\ge 1$, $0<\theta\le 2$ and $1<p<p_\theta$. 
Then there exists $\gamma_2>0$ such that, 
if $\mu$ is a Radon measure in ${\bf R}^N$ satisfying
\begin{equation}
\label{eq:1.9}
\sup_{x\in{\bf R}^N}\mu(B(x,T^{\frac{1}{\theta}}))\le\gamma_2 T^{\frac{N}{\theta}-\frac{1}{p-1}}
\quad\mbox{for some $T>0$}, 
\end{equation}
then problem~\eqref{eq:1.5} possesses a solution in ${\bf R}^N\times[0,T)$.  
\end{theorem}
\begin{theorem}
\label{Theorem:1.4}
Let $N\ge 1$, $0<\theta\le 2$ and $1<\alpha<p$. 
Then there exists $\gamma_3>0$ such that, 
if $\mu$ is a nonnegative measurable function in ${\bf R}^N$ satisfying 
\begin{equation}
\label{eq:1.10}
\sup_{x\in{\bf R}^N}\left[\,\dashint_{B(x,\sigma)}
\mu(y)^\alpha\,dy\,\right]^{\frac{1}{\alpha}}\le\gamma_3\sigma^{-\frac{\theta}{p-1}},
\qquad
0<\sigma<T^{\frac{1}{\theta}},
\end{equation}
for some $T>0$, 
then problem~\eqref{eq:1.5} possesses a solution in ${\bf R}^N\times[0,T)$. 
\end{theorem}
Furthermore, we state the following theorem, 
which is a refinement of Theorem~\ref{Theorem:1.4} in the case $p=p_\theta$ 
and enables us to 
prove the existence of the solution of \eqref{eq:1.5} under assumption~\eqref{eq:1.4}. 
See also Corollary~\ref{Corollary:4.1}. 
\begin{theorem}
\label{Theorem:1.5}
Let $N\ge 1$, $0<\theta\le 2$, $p=p_\theta$ and $\beta>0$. 
For $s>0$, set 
\begin{equation}
\label{eq:1.11}
\Psi_\beta(s):=s[\log (e+s)]^\beta,
\qquad
\rho(s):=
s^{-N}\biggr[\log\biggr(e+\frac{1}{s}\biggr)\biggr]^{-\frac{N}{\theta}}. 
\end{equation}
Then there exists $\gamma_4>0$ such that,  
if $\mu$ is a nonnegative measurable function in ${\bf R}^N$ satisfying 
\begin{equation}
\label{eq:1.12}
\sup_{x\in{\bf R}^N}\Psi_\beta^{-1}\left[\,\dashint_{B(x,\sigma)}
\Psi_\beta(T^\frac{1}{p-1}\mu(y))\,dy\,\right]\le\gamma_4\rho(\sigma T^{-\frac{1}{\theta}}),
\qquad
0<\sigma<T^{\frac{1}{\theta}},
\end{equation}
for some $T>0$, 
then problem~\eqref{eq:1.5} possesses a solution in ${\bf R}^N\times[0,T)$. 
\end{theorem}
\begin{remark}
\label{Remark:1.2}
The assumptions $\alpha>1$ and $\beta>0$ 
are crucial in Theorems~{\rm\ref{Theorem:1.4}} and {\rm\ref{Theorem:1.5}}, respectively. 
Indeed, Theorems~{\rm\ref{Theorem:1.4}} and {\rm\ref{Theorem:1.5}} do not hold 
with $\alpha=1$ and $\beta=0$ in the case $\theta=2$, respectively. 
See {\rm\cite{Takahashi}}.
\end{remark}
As a corollary of Theorems~\ref{Theorem:1.1}, \ref{Theorem:1.4} and \ref{Theorem:1.5}, 
we have
\begin{corollary}
\label{Corollary:1.2}
Let $N\ge 1$, $0<\theta\le 2$ and $p\ge p_\theta$. 
Then there exists $\gamma_*>0$ 
with the following properties: 
\begin{itemize}
\item[{\rm (i)}]
If $p=p_\theta$ and 
$$
\mu(x)=\gamma|x|^{-N}\displaystyle{\biggr[\log\biggr(e+\frac{1}{|x|}\biggr)\biggr]^{-\frac{N}{\theta}-1}}+C
$$
for some $\gamma\ge 0$ and $C\ge 0$, then
\begin{itemize}
  \item problem~\eqref{eq:1.5} possesses a local-in-time solution if $0\le\gamma<\gamma_*$; 
  \item problem~\eqref{eq:1.5} possesses no local-in-time solutions if $\gamma>\gamma_*$. 
\end{itemize}
\item[{\rm (ii)}] 
If $p>p_\theta$ and 
$$
\mu(x)=\gamma|x|^{-\frac{\theta}{p-1}}+C
$$
for some $\gamma\ge 0$ and $C\ge 0$, 
then the same conclusion as in assertion~{\rm (i)} holds. 
Furthermore, if $C=0$ and $\gamma$ is small enough, then 
problem~\eqref{eq:1.5} possesses a global-in-time solution. 
\end{itemize}
\end{corollary}

The rest of this paper is organized as follows. 
In Section~2 we collect some properties of the kernel~$G$ and prove some preliminary lemmas. 
In Section~3 we prove Theorems~\ref{Theorem:1.1} and \ref{Theorem:1.2}. 
In Section~4 we prove Theorems~\ref{Theorem:1.3}, \ref{Theorem:1.4} and \ref{Theorem:1.5}. 
In Section~5, as an application of our theorems, 
we obtain the estimates of the life span of the solution of \eqref{eq:1.1} 
with small initial data. 
\section{Preliminaries}
In this section we collect some properties of the fundamental solution $G$ of \eqref{eq:1.6} 
and recall preliminary lemmas. 
In what follows 
the letter $C$ denotes a generic positive constant depending only on $N$, $\theta$ and $p$.
\vspace{3pt}

Let $N\ge 1$ and $0<\theta\le 2$. 
The fundamental solution $G$ is a positive and smooth function in ${\bf R}^N\times(0,\infty)$ 
and it is represented by 
$$
G(x,t)=
\left\{
\begin{array}{ll}
\displaystyle{(4\pi t)^{-\frac{N}{2}}\exp\left(-\frac{|x|^2}{4t}\right)} & \mbox{if}\quad\theta=2,\vspace{5pt}\\
\displaystyle{\int_0^\infty f_{t,\frac{\theta}{2}}(s)\,(4\pi s)^{-\frac{N}{2}}\exp\left(-\frac{|x|^2}{4s}\right)\,ds} & \mbox{if}\quad 0<\theta<2,
\end{array}
\right. 
$$
where $f_{t,\theta/2}$ is a nonnegative function on $[0,\infty)$ defined by 
$$
f_{t,\frac{\theta}{2}}(s)=\frac{1}{2\pi i}\int_{\sigma-i\infty}^{\sigma+i\infty}
\exp(zs-tz^{\frac{\theta}{2}})\,dz\qquad(\sigma>0,\,\,t>0). 
$$
See \cite[Section~11, Chapter~IX]{Y}. 
Furthermore, $G$ has the following properties, 
\begin{eqnarray}
\label{eq:2.1}
 & & G(x,t)=t^{-\frac{N}{\theta}}G\left(t^{-\frac{1}{\theta}}x,1\right),\\
\label{eq:2.2}
 & & C^{-1}(1+|x|)^{-N-\theta}\le G(x,1)\le C(1+|x|)^{-N-\theta}\quad\mbox{ if $0<\theta<2$},\\
\label{eq:2.3}
 & &  \mbox{$G(\cdot,1)$ is radially symmetric and $G(x,1)\le G(y,1)$ if $|x|\ge |y|$},\\
\label{eq:2.4}
 & & G(x,t)=\int_{{\bf R}^N}G(x-y,t-s)G(y,s)dy,\\
\label{eq:2.5}
 & & \int_{{\bf R}^N}G(x,t)\,dx=1,
\end{eqnarray}
for all $x$, $y\in{\bf R}^N$ and $0<s<t$ 
(see e.g., \cite{BSS, Sugitani}). 

For any locally integrable function $\phi$ in ${\bf R}^N$, 
we often identify $\phi$ with the Radon measure $\phi\,dx$.  
For any Radon measure $\mu$ in ${\bf R}^N$, 
we define
$$
[S(t)\mu](x):=\int_{{\bf R}^N}G(x-y,t)\,d\mu(y),\quad x\in{\bf R}^N,\,\,t>0. 
$$
Then
\begin{equation}
\label{eq:2.6}
\lim_{t\to+0}\|S(t)\eta-\eta\|_{L^\infty({\bf R}^N)}=0,
\qquad\eta\in C_0({\bf R}^N).
\end{equation}
Furthermore, we have 
\begin{lemma}
\label{Lemma:2.1}
There exists a constant $C$ such that
\begin{eqnarray}
\label{eq:2.7}
 & & \|S(t)\mu\|_{L^\infty({\bf R}^N)}\le Ct^{-\frac{N}{\theta}}\sup_{x\in{\bf R}^N}
 \mu(B(x,t^{\frac{1}{\theta}}))
\end{eqnarray}
for any Radon measure $\mu$ in ${\bf R}^N$ and $t>0$. 
\end{lemma}
{\bf Proof.} 
Let $0<\theta<2$, $x \in{\bf R}^N$ and $t>0$. By the Besicovitch covering lemma 
we can find an integer $m$ depending only on $N$ and 
a set $\{x_{k,i}\}_{k=1,\dots,m,\,i\in{\bf N}}\subset{\bf R}^N$ such that 
\begin{equation}
\label{eq:2.8}
B_{k,i}\cap B_{k,j}=\emptyset\quad\mbox{if $i\not=j$}
\qquad\mbox{and}\qquad
{\bf R}^N=\bigcup_{k=1}^m\bigcup_{i=1}^\infty B_{k,i},
\end{equation}
where $B_{k,i}:=B(x_{k,i},t^{1/\theta})$. 
It follows from \eqref{eq:2.1} and \eqref{eq:2.2} that 
\begin{equation}
\label{eq:2.9}
\begin{split}
 [S(t)\mu](x) & \le\sum_{k=1}^m\sum_{i=1}^\infty \int_{B_{k,i}}G(x-y,t)\,d\mu(y) \\
 & \le Ct^{-\frac{N}{\theta}}\sup_{z\in{\bf R}^N}\mu(B(z,t^{\frac{1}{\theta}}))\sum_{k=1}^m \sum_{i=1}^\infty
 \sup_{y\in B_{k,i}}
 \left(1+t^{-\frac{1}{\theta}}|x-y|\right)^{-N-\theta}.
\end{split}
\end{equation}
On the other hand, 
since
\begin{equation*}
\begin{split}
\inf_{y\in B_{k,i}}(1+t^{-\frac{1}{\theta}}|x-y|)
 & \ge 1+\frac{1}{4}t^{-\frac{1}{\theta}}\inf_{y\in B_{k,i}}|x-y|\\
 & \ge 1+\frac{1}{4}t^{-\frac{1}{\theta}}(|x-z|-2t^\frac{1}{\theta})
=\frac{1}{2}+\frac{1}{4}t^{-\frac{1}{\theta}}|x-z|
\end{split}
\end{equation*}
for any $z\in B_{k,i}$, 
we have 
$$
\sup_{y\in B_{k,i}}
\left(1+t^{-\frac{1}{\theta}}|x-y|\right)^{-N-\theta}
\le C\,\,\dashint_{B_{k,i}}
(2+t^{-\frac{1}{\theta}}|x-z|)^{-N-\theta}\,dz. 
$$
This together with \eqref{eq:2.8} and \eqref{eq:2.9} implies that 
\begin{equation*}
\begin{split}
 [S(t)\mu](x) 
 & \le Ct^{-\frac{2N}{\theta}}\sup_{z\in{\bf R}^N}\mu(B(z,t^{\frac{1}{\theta}}))
 \int_{{\bf R}^N}
(2+t^{-\frac{1}{\theta}}|x-z|)^{-N-\theta}\,dz\\
 & \le Ct^{-\frac{N}{\theta}}\sup_{z\in{\bf R}^N}\mu(B(z,t^{\frac{1}{\theta}}))
 \int_{{\bf R}^N}(2+|z|)^{-N-\theta}\,dz
 \le Ct^{-\frac{N}{\theta}}\sup_{z\in{\bf R}^N}\mu(B(z,t^{\frac{1}{\theta}}))
\end{split}
\end{equation*}
for all $x\in{\bf R}^N$ and $t>0$. 
Therefore we obtain \eqref{eq:2.7} in the case $0<\theta<2$. 
Similarly, we have \eqref{eq:2.7} in the case $\theta=2$, 
and the proof is complete. 
$\Box$
\vspace{5pt}

We prove the following three lemmas on Cauchy problem~\eqref{eq:1.5}. 
\begin{lemma}
\label{Lemma:2.2}
Let $\mu$ be a Radon measure in ${\bf R}^N$ and $0<T\le\infty$. 
Assume that there exists a supersolution $v$ of \eqref{eq:1.5} in ${\bf R}^N\times[0,T)$.  
Then there exists a minimal solution of  \eqref{eq:1.5} in ${\bf R}^N\times[0,T)$. 
\end{lemma}
{\bf Proof.} 
Set $u_1:=S(t)\mu$. 
Define $u_n$ $(n=2,3,\dots)$ inductively by 
\begin{equation}
\label{eq:2.10}
u_n(t):=S(t)\mu+\int_0^t S(t-s)u_{n-1}(s)^p\,ds.
\end{equation}
Let $v$ be a supersolution of \eqref{eq:1.5} in ${\bf R}^N\times[0,T)$, 
where $0<T\le\infty$. 
Then it follows inductively that 
$$
0\le u_1(x,t)\le u_2(x,t)\le\cdots\le u_n(x,t)\le\cdots\le v(x,t)<\infty
$$
for almost all $x\in{\bf R}^N$ and $t\in(0,T)$. 
This means that 
$$
u(x,t):=\lim_{n\to\infty}u_n(x,t)\le v(x,t)
$$
for almost all $x\in{\bf R}^N$ and $t\in(0,T)$. 
Furthermore, by \eqref{eq:2.10} 
we see that $u$ satisfies \eqref{eq:1.7} 
for almost all $x\in{\bf R}^N$ and $t\in(0,T)$. 
In addition, we easily see that $u$ is a minimal solution of \eqref{eq:1.5} in ${\bf R}^N\times[0,T)$. 
Thus Lemma~\ref{Lemma:2.2} follows. 
$\Box$
\begin{lemma}
\label{Lemma:2.3}
Let $u$ be a solution of \eqref{eq:1.1} in ${\bf R}^N\times(0,T)$, where $0<T<\infty$. 
Then 
\begin{equation}
\label{eq:2.11}
\underset{0<t<T-\epsilon}{\mbox{{\rm ess sup}}}\,\int_{B(0,R)}u(y,t)\,dy<\infty
\end{equation}
for all $R>0$ and $0<\epsilon<T$. 
Furthermore, 
there exists a unique Radon measure $\mu$ in ${\bf R}^N$ such that
\begin{equation}
\label{eq:2.12}
\underset{t\to +0}{\mbox{{\rm ess lim}}}
\int_{{\bf R}^N}u(y,t)\eta(y)\,dy=\int_{{\bf R}^N}\eta(y)\,d\mu(y)
\end{equation}
for all $\eta\in C_0({\bf R}^N)$. 
\end{lemma}
{\bf Proof.}
Let $0<\epsilon<T/2$ and $R>0$. 
It follows from Definition~\ref{Definition:1.1}~(i) that 
\begin{equation*}
\begin{split}
\infty>u(x,t) & \ge\int_{B(x,R)}G(x-y,t-\tau)u(y,\tau)\,dy\\
 & \ge\inf_{z\in B(0,R), \epsilon/2<s<T}G(z,s)
\int_{B(x,R)}u(y,\tau)\,dy
\ge C_\epsilon\int_{B(x,R)}u(y,\tau)\,dy,
\end{split}
\end{equation*}
for almost all $x\in{\bf R}^N$, $\tau\in(0,T-\epsilon)$ and $t\in(T-\epsilon/2,T)$, 
where $C_\epsilon$ is a positive constant independent of $x$ and $\tau$. 
This implies \eqref{eq:2.11}. 
Then, 
applying the weak compactness of Radon measures (see e.g., \cite[Section~1.9]{EG}), 
we can find a sequence $\{t_j\}$ with $\lim_{j\to\infty}t_j=0$ 
and a Radon measure $\mu$ in ${\bf R}^N$ such that 
\begin{equation}
\label{eq:2.13}
\lim_{j\to\infty}\int_{{\bf R}^N}u(y,t_j)\eta(y)\,dy=\int_{{\bf R}^N}\eta(y)\,d\mu(y)
\end{equation}
for all $\eta\in C_0({\bf R}^N)$. 

We prove the uniqueness of the Radon measures satisfying \eqref{eq:2.13}. 
Assume that 
there exist a sequence $\{s_j\}$ with $\lim_{j\to\infty}s_j=0$ and a Radon measure $\mu'$ in ${\bf R}^N$ such that 
\begin{equation}
\label{eq:2.14}
\lim_{j\to\infty}\int_{{\bf R}^N}u(y,s_j)\eta(y)\,dy=\int_{{\bf R}^N}\eta(y)\,d\mu'(y)
\end{equation}
for all $\eta\in C_0({\bf R}^N)$. 
Let $\{s_{j'}\}$ be a subsequence of $\{s_j\}$ such that $t_j>s_{j'}$ for $j=1,2,\dots$. 
By using Definition~\ref{Definition:1.1}~(i) again 
we see that 
\begin{equation}
\label{eq:2.15}
u(x,t_j)\ge\int_{{\bf R}^N} G(x-y,t_j-s_{j'})u(y,s_{j'})\,dy,
\qquad j=1,2,\dots.
\end{equation}
For any $R>0$,  
let $\zeta\in C_0(B(0,R))$ be such that $\zeta\ge 0$ in $B(0,R)$. 
Set $\zeta(x,t):=[S(t)\zeta](x)$. 
By \eqref{eq:2.3} and \eqref{eq:2.15} 
we have 
\begin{equation*}
\begin{split}
 & \int_{{\bf R}^N} u(x,t_j)\zeta(x)\,dx
 \ge\int_{{\bf R}^N}\left(\int_{{\bf R}^N} G(x-y,t_j-s_{j'})\zeta(x)\,dx\right)\,u(y,s_{j'})\,dy\\
 & =\int_{{\bf R}^N} \zeta(y,t_j-s_{j'})u(y,s_{j'})\,dy
 \ge\int_{B(0,R)} \zeta(y,t_j-s_{j'})u(y,s_{j'})\,dy\\
 & \ge \int_{{\bf R}^N}\zeta(y)u(y,s_{j'})\,dy-\|\zeta(t_j-s_{j'})-\zeta\|_{L^\infty(B(0,R))}
 \int_{B(0,R)}u(y,s_{j'})\,dy.
\end{split} 
\end{equation*}
Then, letting $j\to\infty$, 
by \eqref{eq:2.6}, \eqref{eq:2.11}, \eqref{eq:2.13} and \eqref{eq:2.14} 
we obtain 
$$
\int_{{\bf R}^N}\zeta\,d\mu\ge\int_{{\bf R}^N}\zeta\,d\mu'.
$$
Similarly, it follows that 
$$
\int_{{\bf R}^N}\zeta\,d\mu'\ge\int_{{\bf R}^N}\zeta\,d\mu.
$$
Since $\zeta$ is arbitrary, we deduce that $\mu=\mu'$ in ${\bf R}^N$. 
Therefore we obtain the uniqueness of the Radon measures satisfying \eqref{eq:2.13}. 
Then \eqref{eq:2.12} follows from \eqref{eq:2.13}. Thus the proof is complete. 
$\Box$
\begin{lemma}
\label{Lemma:2.4}
Let $\mu$ be a Radon measure in ${\bf R}^N$. 
Let $u$ be a solution of \eqref{eq:1.1} in ${\bf R}^N\times[0,T)$ 
with $u(\cdot,0)=\mu$, where $0<T<\infty$. Then  
\begin{equation}
\label{eq:2.16}
\underset{t\to+0}{\mbox{{\rm ess lim}}}
\int_{{\bf R}^N}u(x,t)\eta(x)\,dx=\int_{{\bf R}^N}\eta(x)\,d\mu(x)
\end{equation}
for all $\eta\in C_0({\bf R}^N)$. 
\end{lemma}
{\bf Proof.}
It suffices to prove \eqref{eq:2.16} 
for all $\eta\in C_0({\bf R}^N)$ with $\eta\ge 0$ in ${\bf R}^N$. 
Let $0<\theta<2$ and $R\ge 1$ be such that $\mbox{supp}\,\eta\subset B(0,R)$. 
In the proof, 
the letter $C$ denotes a generic positive constant depending only on $N$, $\theta$, $p$, $R$ and $T$. 

It follows from the Fubini theorem and \eqref{eq:2.4} that 
$u$ is a solution of \eqref{eq:1.1} in ${\bf R}^N\times(0,T)$ in the sense of Definition~\ref{Definition:1.1}~(i). 
Then, by Lemma~\ref{Lemma:2.3} 
we see that 
\begin{equation}
\label{eq:2.17}
\underset{0<t<T-\epsilon}{\mbox{{\rm ess sup}}}\,\int_{B(0,R)}u(y,t)\,dy<\infty
\quad\mbox{for all $R>0$ and $0<\epsilon<T$}
\end{equation}
and we can find a unique Radon measure $\mu'$ in ${\bf R}^N$ such that 
\begin{equation}
\label{eq:2.18}
\underset{t\to +0}{\mbox{{\rm ess lim}}}
\int_{{\bf R}^N}u(y,t)\eta(y)\,dy=\int_{{\bf R}^N}\eta(y)\,d\mu'(y).
\end{equation}
Furthermore, 
\begin{equation}
\label{eq:2.19}
\begin{split}
\int_{{\bf R}^N}u(x,t)\eta(x)\,dx
 & =\int_{{\bf R}^N}\int_{{\bf R}^N}G(x-y,t-\tau)u(y,\tau)\eta(x)\,dx\,dy\\
 & \qquad
 +\int_\tau^t\int_{{\bf R}^N}\int_{{\bf R}^N}G(x-y,t-s)u(y,s)^p\eta(x)\,dx\,dy\,ds\\
 & \ge\int_{B(0,R)}\int_{{\bf R}^N}G(x-y,t-\tau)u(y,\tau)\eta(x)\,dx\,dy\\
  & \qquad
 +\int_\tau^t\int_{{\bf R}^N}\int_{{\bf R}^N}G(x-y,t-s)u(y,s)^p\eta(x)\,dx\,dy\,ds
\end{split}
\end{equation}
for almost all $0<\tau<t<T/2$. 
Set $\eta(x,t):=[S(t)\eta](x)$. 
It follows from \eqref{eq:2.3} that 
\begin{equation*}
\begin{split}
 & \int_{B(0,R)}\int_{{\bf R}^N}G(x-y,t-\tau)u(y,\tau)\eta(x)\,dx\,dy
=\int_{B(0,R)}\eta(y,t-\tau)u(y,\tau)\,dy\\
 & \ge\int_{B(0,R)}\eta(y)u(y,\tau)\,dy-\|\eta(t-\tau)-\eta\|_{L^\infty(B(0,R))}\int_{B(0,R)}u(y,\tau)\,dy,
\end{split}
\end{equation*}
which together with \eqref{eq:2.17} and \eqref{eq:2.18} implies that 
\begin{equation}
\label{eq:2.20}
\begin{split}
 & \underset{\tau\to+0}{\mbox{{\rm ess liminf}}}\int_{B(0,R)}\int_{{\bf R}^N}G(x-y,t-\tau)u(y,\tau)\eta(x)\,dx\,dy\\
 & \qquad\qquad
\ge\int_{{\bf R}^N}\eta(y)\,d\mu'(y)-C\|\eta(t)-\eta\|_{L^\infty(B(0,R))}
\end{split}
\end{equation}
for almost all $0<t<T/2$. 
By \eqref{eq:2.19} and \eqref{eq:2.20} we obtain 
\begin{equation*}
\begin{split}
\int_{{\bf R}^N}u(x,t)\eta(x)\,dx
 & \ge\int_{{\bf R}^N}\eta(y)\,d\mu'(y)-C\|\eta(t)-\eta\|_{L^\infty(B(0,R))}\\
 & \qquad
 +\int_0^t\int_{{\bf R}^N}\int_{{\bf R}^N}G(x-y,t-s)u(y,s)^p\eta(x)\,dx\,dy\,ds
\end{split}
\end{equation*}
for almost all $0<t<T/2$. 
This together with \eqref{eq:2.6} and \eqref{eq:2.18} implies that 
\begin{equation}
\label{eq:2.21}
\underset{t\to+0}{\mbox{{\rm ess lim}}}
\int_0^t\int_{{\bf R}^N}\int_{{\bf R}^N}G(x-y,t-s)u(y,s)^p\eta(x)\,dx\,dy\,ds=0. 
\end{equation}

On the other hand, 
it follows from \eqref{eq:2.1} and \eqref{eq:2.2} that
\begin{equation}
\label{eq:2.22}
\eta(y,t) 
\le C t\int_{B(0,R)} (t^\frac{1}{\theta} + |x-y|)^{-N-\theta} \eta(x) \,dx
\le C\|\eta\|_{L^\infty({\bf R}^N)}(1+|y|)^{-N-\theta}
\end{equation}
for $y\in{\bf R}^N$ with $|y|\ge 2R\ge 2$ and $0<t\le T/2$. 
By \eqref{eq:2.5} we see that $\|\eta(t)\|_{L^\infty({\bf R}^N)}\le\|\eta\|_{L^\infty({\bf R}^N)}$ for $t>0$. 
This together with \eqref{eq:2.22} implies that 
\begin{equation}
\label{eq:2.23}
0\le \eta(y,t) 
\le C\|\eta\|_{L^\infty({\bf R}^N)}(1+|y|)^{-N-\theta}
\end{equation}
for $y\in{\bf R}^N$ and $0<t\le T/2$. 
It follows from \eqref{eq:2.1} and \eqref{eq:2.2} that 
\begin{equation}
\label{eq:2.24}
\int_{B(0,1)} G(x-y,\tau)\, dx 
\ge C\int_{B(0,1)} (1+|x-y|)^{-N-\theta}\, dx\\
\ge C(1+|y|)^{-N-\theta}
\end{equation}
for $y\in{\bf R}^N$ and $\tau\in(T/4,T/2)$. 
By \eqref{eq:1.7}, \eqref{eq:2.17}, \eqref{eq:2.23} and \eqref{eq:2.24} we see that 
\begin{equation}
\label{eq:2.25}
\begin{split}
\infty> \underset{T/4<\tau<T/2}{\mbox{{\rm ess sup}}}\int_{B(0,1)} u(x,\tau)\,dx 
 & \ge \underset{T/4<\tau<T/2}{\mbox{{\rm ess sup}}}\,\int_{B(0,1)} \int_{{\bf R}^N} G(x-y,\tau)\,d\mu(y)\,dx\\
 & \ge C\int_{{\bf R}^N} (1+|y|)^{-N-\theta}\,d\mu(y).
\end{split} 
\end{equation}
By \eqref{eq:2.6}, \eqref{eq:2.23} and \eqref{eq:2.25} 
we apply the Lebesgue dominated convergence theorem to obtain  
$$
\lim_{t\to +0}
\int_{{\bf R}^N}(\eta(y,t)-\eta(y))\,d\mu(y)=0,
$$
which implies that  
\begin{equation}
\label{eq:2.26}
\begin{split}
 & \underset{t\to+0}{\mbox{{\rm ess lim}}}
 \int_{{\bf R}^N}\int_{{\bf R}^N}G(x-y,t)\eta(x)\,dx\,d\mu(y)\\
 & =\int_{{\bf R}^N}\eta(y)\,d\mu(y)
 +\lim_{t\to +0}\int_{{\bf R}^N}(\eta(y,t)-\eta(y))\,d\mu(y)
 =\int_{{\bf R}^N}\eta(y)\,d\mu(y).
\end{split}
\end{equation}
Therefore we deduce from \eqref{eq:1.7}, \eqref{eq:2.19}, \eqref{eq:2.21} and \eqref{eq:2.26} that 
$$
\underset{t\to+0}{\mbox{{\rm ess lim}}}\int_{{\bf R}^N} u(x,t)\eta(x)\,dx=\int_{{\bf R}^N}\eta(y)\,d\mu(y). 
$$
Similarly, we have \eqref{eq:2.16} in the case $\theta=2$,
and Lemma~\ref{Lemma:2.4} follows. 
$\Box$
\section{Proof of Theorems~\ref{Theorem:1.1} and \ref{Theorem:1.2}}
We first prove the following lemma by using the argument in \cite[Theorem~5]{W1}. 
\begin{lemma}
\label{Lemma:3.1}
Let $u$ be a solution of \eqref{eq:1.1} in ${\bf R}^N\times(0,T)$, where $0<T<\infty$. 
Then there exists a constant $\kappa$ depending only on $p$ such that 
\begin{equation}
\label{eq:3.1}
\|S(t)u(\tau)\|_{L^\infty({\bf R}^N)}\le\kappa t^{-\frac{1}{p-1}}
\end{equation}
for almost all $t>0$ and $\tau>0$ with $t+\tau<T$. 
\end{lemma}
{\bf Proof.}
It follows from Definition~\ref{Definition:1.1}~(i) that 
\begin{equation}
\label{eq:3.2}
u(x,t+\tau)\ge [S(t)u(\tau)](x)
\end{equation}
for almost all $x\in{\bf R}^N$, $t>0$ and $\tau>0$ with $t+\tau<T$. 
By the Jensen inequality, Definition~\ref{Definition:1.1}~(i), \eqref{eq:2.5} and \eqref{eq:3.2}
we see that 
\begin{equation}
\label{eq:3.3}
\begin{split}
\infty>u(x,t+\tau) & \ge \int_\tau^{t+\tau}\int_{{\bf R}^N}G(x-y,t+\tau-s)[S(s-\tau)u(\tau)](y)^p\,dy\,ds\\
 & \ge\int_\tau^{t+\tau}\biggr(\int_{{\bf R}^N}G(x-y,t+\tau-s)[S(s-\tau)u(\tau)](y)\,dy\biggr)^p\,ds\\
 & =\int_\tau^{t+\tau}[S(t)u(\tau)](x)^p\,ds
=t[S(t)u(\tau)](x)^p
\end{split}
\end{equation}
for almost all $x\in{\bf R}^N$, $t>0$ and $\tau>0$ with $t+\tau<T$. 
Similarly, 
by using \eqref{eq:3.3}, instead of \eqref{eq:3.2}, we obtain 
\begin{equation*}
\begin{split}
\infty>u(x,t+\tau) & \ge \int_\tau^{t+\tau}\int_{{\bf R}^N}G(x-y,t+\tau-s)(s-\tau)^p[S(s-\tau)u(\tau)](y)^{p^2}\,dy\,ds\\
 & \ge\int_\tau^{t+\tau} (s-\tau)^p\biggr(\int_{{\bf R}^N}G(x-y,t+\tau-s)[S(s-\tau)u(\tau)](y)\,dy\biggr)^{p^2}\,ds\\
 & =\int_\tau^{t+\tau}(s-\tau)^p[S(t)u(\tau)](x)^{p^2}\,ds
=\frac{1}{p+1}t^{p+1}[S(t)u(\tau)](x)^{p^2}
\end{split}
\end{equation*}
for almost all $x\in{\bf R}^N$, $t>0$ and $\tau>0$ with $t+\tau<T$. 
Then, by the same argument as in the proof of \cite[Theorem~5]{W1} 
we can find a constant $\kappa$ depending only on $p$  such that 
$$
t^{\frac{1}{p-1}}[S(t)u(\tau)](x)\le\kappa
$$
for almost all $x\in{\bf R}^N$, $t>0$ and $\tau>0$ with $t+\tau<T$.  
This implies \eqref{eq:3.1}, and 
the proof is complete. 
$\Box$\vspace{5pt}

Next we refine \eqref{eq:2.11} and obtain the following lemma. 
\begin{lemma}
\label{Lemma:3.2}
Let $u$ be a solution of \eqref{eq:1.1} in ${\bf R}^N\times(0,T)$, 
where $0<T<\infty$. 
Let $z\in{\bf R}^N$ and $\rho>0$ with $(2\rho)^\theta<T$. 
Then there exists a constant $c_*$ depending only on $N$ and $\theta$ such that 
\begin{equation}
\label{eq:3.4}
u(x+z,\tau+(2\rho)^\theta)\ge c_*G(x,\rho^\theta)\int_{B(z,\rho)}u(y,\tau)\,dy
\end{equation}
for almost all $x\in{\bf R}^N$ and $\tau\in(0,T- (2\rho)^\theta)$. 
\end{lemma}
{\bf Proof.}
Let $z\in{\bf R}^N$. 
Since 
$$
|x+z-y|\le 2|x|\quad\mbox{if}\quad |x|>\rho\quad\mbox{and}\quad|y-z|<\rho,
$$
by Definition~\ref{Definition:1.1}~(i),  \eqref{eq:2.1} and \eqref{eq:2.3} 
we obtain 
\begin{equation}
\label{eq:3.5}
\begin{split}
 & u(x+z,\tau+(2\rho)^\theta)
\ge\int_{B(z,\rho)}G(x+z-y,(2\rho)^\theta)u(y,\tau)\,dy\\
 & \qquad\quad
\ge \int_{B(z,\rho)}G(2x,(2\rho)^\theta)u(y,\tau)\,dy=2^{-N}G(x,\rho^\theta)\int_{B(z,\rho)}u(y,\tau)\,dy
\end{split}
\end{equation}
for almost all $x\in{\bf R}^N$ with $|x|>\rho$ and $\tau\in(0,T-(2\rho)^\theta)$. 
Similarly, we have  
\begin{equation}
\label{eq:3.6}
\begin{split}
u(x+z,\tau+(2\rho)^\theta) & \ge\int_{B(z,\rho)}G(x+z-y,(2\rho)^\theta)u(y,\tau)\,dy\\
 & \ge \min_{|z|\le 2\rho}G(z,(2\rho)^\theta)\int_{B(z,\rho)}u(y,\tau)\,dy\\
 & \ge \frac{\min_{|z|\le 2\rho}G(z,(2\rho)^\theta)}{\max_{|z|\le \rho}G(z,\rho^\theta)}
G(x,\rho^\theta)\int_{B(z,\rho)}u(y,\tau)\,dy\\
 & =2^{-N}\frac{\min_{|z|\le 1}G(z,1)}{\max_{|z|\le 1}G(z,1)}G(x,\rho^\theta)\int_{B(z,\rho)}u(y,\tau)\,dy
\end{split}
\end{equation}
for almost all $x\in{\bf R}^N$ with $|x|\le\rho$ and $\tau\in(0,T-(2\rho)^\theta)$. 
By \eqref{eq:3.5} and \eqref{eq:3.6} we obtain \eqref{eq:3.4}, 
and the proof is complete. 
$\Box$\vspace{5pt}

Now we are ready to prove Theorem~\ref{Theorem:1.1} in the case $p\not=p_\theta$.
\vspace{5pt}
\newline
{\bf Proof of Theorem~\ref{Theorem:1.1} in the case $p\not=p_\theta$.}
Let $u$ be a solution of \eqref{eq:1.1} in ${\bf R}^N\times(0,T)$, where $0<T<\infty$. 
Let $\rho>0$ be such that $2(2\rho)^\theta<T$. 
For any $z\in{\bf R}^N$, set 
$$
\tilde{u}(x,t):=u(x+z,t+(2\rho)^\theta)
$$
for almost all $x\in{\bf R}^N$ and $t\in(0,T-(2\rho)^\theta)$. 
Since $\tilde{u}$ is a solution of \eqref{eq:1.1} in ${\bf R}^N\times(0,T-(2\rho)^\theta)$, 
Lemma~\ref{Lemma:3.1} implies 
\begin{equation}
\label{eq:3.7}
\|S((2\rho)^\theta)\tilde{u}(\tau))\|_{L^\infty({\bf R}^N)}\le\kappa (2\rho)^{-\frac{\theta}{p-1}}
\end{equation}
for almost all $\rho>0$ with $2(2\rho)^\theta<T$ and $\tau>0$ with $\tau+(2\rho)^\theta<T$. 
Furthermore, 
it follows from Lemma~\ref{Lemma:3.2} and \eqref{eq:2.4} that  
\begin{equation*}
\begin{split}
 & [S((2\rho)^\theta)\tilde{u}(\tau)](x)
=\int_{{\bf R}^N}G(x-y,(2\rho)^\theta)u(y+z,\tau+(2\rho)^\theta)\,dy\\
 & \qquad\quad
\ge C\int_{{\bf R}^N}G(x-y,(2\rho)^\theta)G(y,\rho^\theta)\,dy \int_{B(z,\rho)}u(y,\tau)\,dy\\
 & \qquad\quad
=CG(x,(2\rho)^\theta+\rho^\theta)\int_{B(z,\rho)}u(y,\tau)\,dy
\end{split}
\end{equation*}
for all $\rho>0$ with $2(2\rho)^\theta<T$ and almost all $\tau>0$ with $\tau+(2\rho)^\theta<T$. 
This together with \eqref{eq:2.1} and \eqref{eq:2.3} implies that 
\begin{equation}
\label{eq:3.8}
\begin{split}
 & \|S((2\rho)^\theta)\tilde{u}(\tau))\|_{L^\infty({\bf R}^N)}
\ge CG(0,(2\rho)^\theta+\rho^\theta)\int_{B(z,\rho)}u(y,\tau)\,dy\\
 & \qquad\quad
\ge C[(2\rho)^\theta+\rho^\theta)]^{-\frac{N}{\theta}}G(0,1)
\int_{B(z,\rho)}u(y,\tau)\,dy\\
 & \qquad\quad
\ge C\rho^{-N}\int_{B(z,\rho)}u(y,\tau)\,dy
\end{split}
\end{equation}
for all $\rho>0$ with $2(2\rho)^\theta<T$ and almost all $\tau>0$ with $\tau+(2\rho)^\theta<T$. 
We deduce from  \eqref{eq:3.7} and \eqref{eq:3.8} that 
\begin{equation}
\label{eq:3.9}
\int_{B(z,\rho)}u(y,\tau)\,dy\le C\rho^{N-\frac{\theta}{p-1}}
\end{equation}
for all $z\in{\bf R}^N$ and $\rho>0$ with $2(2\rho)^\theta<T$ and 
for almost all $\tau>0$ with $\tau+(2\rho)^\theta<T$. 

On the other hand, by Lemma~\ref{Lemma:2.3} 
we can find a unique Radon measure $\mu$ such that 
\begin{equation}
\label{eq:3.10}
\underset{t\to+0}{\mbox{{\rm ess lim}}}
\int_{{\bf R}^N}u(y,t)\eta(y)\,dy=\int_{{\bf R}^N}\eta(y)\,d\mu(y),
\qquad
\eta\in C_0({\bf R}^N). 
\end{equation}
Let $\zeta\in C_0({\bf R}^N)$ be such that 
\begin{equation}
\label{eq:3.11}
\zeta=1\quad\mbox{in}\quad B(0,\rho/2),
\quad
0\le\zeta\le 1\quad\mbox{in}\quad {\bf R}^N,
\quad
\zeta=0\quad\mbox{outside}\quad B(0,\rho).
\end{equation}
By \eqref{eq:3.9}, \eqref{eq:3.10} and \eqref{eq:3.11} 
we see that 
\begin{equation}
\label{eq:3.12}
C\rho^{N-\frac{\theta}{p-1}}\ge\underset{\tau\to+0}{\mbox{{\rm ess lim}}}\int_{B(z,\rho)}u(y,\tau)\zeta(y-z)\,dy
=\int_{{\bf R}^N}\zeta(y-z)\,d\mu(y)\ge\mu(B(z,\rho/2))
\end{equation}
for all $z\in{\bf R}^N$ and $\rho>0$ with $2(2\rho)^\theta<T$. 
Setting $\sigma:=2^{(1+\theta)/\theta}\rho$, 
we obtain 
$$
\sup_{z\in{\bf R}^N}\mu\left(B(z,2^{-\frac{1+2\theta}{\theta}}\sigma)\right)\le C\sigma^{N-\frac{\theta}{p-1}}
$$
for all $0<\sigma<T^{1/\theta}$. 
Then we can find an integer $m$ such that 
\begin{equation}
\label{eq:3.13}
\sup_{z\in{\bf R}^N}\mu(B(z,\sigma))
\le m\sup_{z\in{\bf R}^N}\mu\left(B(z,2^{-\frac{1 + \theta}{\theta}}\sigma)\right)
\le Cm\sigma^{N-\frac{\theta}{p-1}}
\end{equation}
for all $0<\sigma<T^{1/\theta}$ (see \cite[Lemma~2.1]{IS}). 
This implies that 
\begin{eqnarray}
\notag
 & & \sup_{z\in{\bf R}^N}\mu(B(z,T^{\frac{1}{\theta}}))
\le CmT^{\frac{N}{\theta}-\frac{1}{p-1}}\quad\mbox{if $1<p<p_\theta$},\\
\label{eq:3.14}
 & & \sup_{z\in{\bf R}^N}\mu(B(z,\sigma))
 \le Cm\sigma^{N-\frac{\theta}{p-1}}
 \quad\mbox{for all $0<\sigma<T^{\frac{1}{\theta}}$ if $p\ge p_\theta$}.
\end{eqnarray}
Therefore we obtain the desired result for $p\not=p_\theta$, 
and the proof of Theorem~\ref{Theorem:1.1} in the case $p\not=p_\theta$ is complete. 
$\Box$\vspace{5pt}

We improve \eqref{eq:3.14} in the case $p=p_\theta$ 
and complete the proof of Theorem~\ref{Theorem:1.1}. 
For this aim, we prove the following lemma. 
\begin{lemma}
\label{Lemma:3.3}
Let $p=p_\theta$. 
Let $u$ be a solution of \eqref{eq:1.1} in ${\bf R}^N\times(0,T)$, 
where $0<T<\infty$.
Then there exist positive constants $C$ and $\nu$ depending only on $N$ and $\theta$ 
such that 
$$
\sup_{z\in{\bf R}^N}\int_{B(z,\rho)}u(y,\tau)\,dy
\le C\log\left(e+\frac{T^{\frac{1}{\theta}}}{\rho}\right)^{-\frac{N}{\theta}}
$$
for all $\rho>0$ with $0<\rho^\theta<\nu T$ 
and for almost all $\tau\in(0,\rho^\theta)$. 
\end{lemma}
{\bf Proof.}
Let $\nu$ be a sufficiently small constant. 
Let  
\begin{equation}
\label{eq:3.15}
0<\rho<(\nu T)^{\frac{1}{\theta}}. 
\end{equation}
For any $z\in{\bf R}^N$, we set 
$v(x,t):=u(x+z,t+(2\rho)^\theta)$ 
for almost all $x\in{\bf R}^N$ and $t\in(0,T-(2\rho)^\theta)$. 
Then it follows from Definition~\ref{Definition:1.1}~(i) that 
\begin{equation}
\label{eq:3.16}
v(x,t)=\int_{{\bf R}^N}G(x-y,t-\tau)v(y,\tau)\,dy
+\int_\tau^t\int_{{\bf R}^N}G(x-y,t-s)v(y,s)^p\,dy\,ds<\infty
\end{equation}
for almost all $x\in{\bf R}^N$ and $0<\tau<t<T-(2\rho)^\theta$. 
On the other hand, in the case $0<\theta<2$, 
by \eqref{eq:2.1} and \eqref{eq:2.2} we have 
$$
G(x-y,\tau)\ge C\tau^{-\frac{N}{\theta}}\left(1+\frac{|x|+|y|}{\tau^{\frac{1}{\theta}}}\right)^{-N-\theta}
\ge C\tau^{-\frac{N}{\theta}}\left(2+\frac{|y|}{\tau^{\frac{1}{\theta}}}\right)^{-N-\theta}
\ge CG(y,\tau)
$$
for all $x\in{\bf R}^N$ with $|x|<\tau^{1/\theta}$, $y\in{\bf R}^N$ and $\tau>0$. 
Then, applying \eqref{eq:3.16} with $t=2\tau$, 
we see that
\begin{equation}
\label{eq:3.17}
\int_{{\bf R}^N}G(y,\tau)v(y,\tau)\,dy<\infty
\end{equation}
for almost all $\tau\in(0,[T-(2\rho)^\theta]/2)$ in the case $0<\theta<2$. 
In the case $\theta=2$, since 
$$
G(x-y,2\tau)\ge (8\pi\tau)^{-\frac{N}{2}}
\exp\left(-\frac{2|x|^2+2|y|^2}{8\tau}\right)
\ge C(4\pi\tau)^{-\frac{N}{2}}\exp\left(-\frac{|y|^2}{4\tau}\right)
$$
for all $x\in{\bf R}^N$ with $|x|<\tau^{1/2}$, $y\in{\bf R}^N$ and $\tau>0$, 
applying \eqref{eq:3.16} with $t=3\tau$,
we have 
\begin{equation}
\label{eq:3.18}
\int_{{\bf R}^N}G(y,\tau)v(y,\tau)\,dy<\infty
\end{equation}
for almost all $\tau\in(0,[T-(2\rho)^\theta]/3)$.
Furthermore, 
by Lemma~\ref{Lemma:3.2}, \eqref{eq:2.4} and \eqref{eq:3.16} 
we have 
\begin{equation}
\label{eq:3.19}
\begin{split}
 & v(x,t)-\int_\tau^t\int_{{\bf R}^N}G(x-y,t-s)v(y,s)^p\,dy\,ds\\
 & =\int_{{\bf R}^N}G(x-y,t-\tau)u(y+z,\tau+(2\rho)^\theta)\,dy\\
 & \ge c_*\int_{B(z,\rho)}u(y,\tau)\,dy
\int_{{\bf R}^N}G(x-y,t-\tau)G(y,\rho^\theta)\,dy\\
 & =c_*\int_{B(z,\rho)}u(y,\tau)\,dy\,G(x,t-\tau+\rho^\theta)
\end{split}
\end{equation}
for almost all $x\in{\bf R}^N$ and $0<\tau<t<T-(2\rho)^\theta$, 
where $c_*$ is the constant as in Lemma~\ref{Lemma:3.2}. 
Set 
$$
w(t):=\int_{{\bf R}^N}v(x,t)G(x,t)\,dx,\qquad 
M_\tau:=\int_{B(z,\rho)}u(y,\tau)\,dy.
$$
Then it follows from \eqref{eq:2.4}, \eqref{eq:3.17}, \eqref{eq:3.18} and \eqref{eq:3.19} that 
\begin{equation}
\label{eq:3.20}
\begin{split}
\infty>w(t) & \ge c_*M_\tau\int_{{\bf R}^N}G(x,t-\tau+\rho^\theta)G(x,t)\,dx\\
 & \qquad\quad
 +\int_{{\bf R}^N}\int_\tau^t\int_{{\bf R}^N}G(x-y,t-s)G(x,t)v(y,s)^p\,dy\,ds\,dx\\
 & \ge c_*M_\tau G(0,2t-\tau+\rho^\theta)+
\int_{\rho^\theta}^t\int_{{\bf R}^N}G(y,2t-s)v(y,s)^p\,dy\,ds
\end{split}
\end{equation}
for almost all $0<\tau<\rho^\theta$ and $\rho^\theta<t<[T-(2\rho)^\theta]/3$. 
On the other hand, 
it follows that from \eqref{eq:2.1} and \eqref{eq:2.3} that
\begin{equation}
\label{eq:3.21}
\begin{split}
G(y,2t-s) & =(2t-s)^{-\frac{N}{\theta}}G\left(\frac{y}{(2t-s)^{\frac{1}{\theta}}},1\right)\\
 & \ge\left(\frac{s}{2t}\right)^{\frac{N}{\theta}}s^{-\frac{N}{\theta}}G\left(\frac{y}{s^{\frac{1}{\theta}}},1\right)
=\left(\frac{s}{2t}\right)^{\frac{N}{\theta}}G(y,s)
\end{split}
\end{equation}
for $y\in{\bf R}^N$ and $0<s<t$. 
Applying the Jensen inequality, 
by \eqref{eq:2.1}, \eqref{eq:2.5}, \eqref{eq:3.20} and \eqref{eq:3.21} 
we obtain 
\begin{equation}
\label{eq:3.22}
\begin{split}
\infty & >w(t)\ge c_*M_\tau G(0,2t-\tau+\rho^\theta)+
\int_{\rho^\theta}^t\left(\frac{s}{2t}\right)^{\frac{N}{\theta}}\int_{{\bf R}^N}G(y,s)v(y,s)^p\,dy\,ds\\
 & \ge c_*M_\tau(2t-\tau+\rho^\theta)^{-\frac{N}{\theta}}G(0,1)+
\int_{\rho^\theta}^t\left(\frac{s}{2t}\right)^{\frac{N}{\theta}}\left(\int_{{\bf R}^N}G(y,s)v(y,s)\,dy\right)^p\,ds\\
& \ge c_*M_\tau 3^{-\frac{N}{\theta}}G(0,1)t^{-\frac{N}{\theta}}
+2^{-\frac{N}{\theta}}t^{-\frac{N}{\theta}}\int_{\rho^\theta}^t
s^{\frac{N}{\theta}}w(s)^p\,ds
\end{split}
\end{equation}
for almost all $0<\tau<\rho^\theta$ and $\rho^\theta<t<[T-(2\rho)^\theta]/3$. 

For $k=1,2,\dots$, we define $\{a_k\}$ inductively by 
\begin{equation}
\label{eq:3.23}
a_1:=c_*3^{-\frac{N}{\theta}}G(0,1),\qquad
a_{k+1}:=2^{- \frac{N}{\theta}}a_k^p\,\frac{p-1}{p^k-1}
\quad(k=1,2,\dots).  
\end{equation}
Furthermore, set 
\begin{equation}
\label{eq:3.24}
f_k(t):=a_kM_\tau^{p^{k-1}}t^{-\frac{N}{\theta}}\left(\log\frac{t}{\rho^\theta}\right)^{\frac{p^{k-1}-1}{p-1}},
\qquad k=1,2,\dots.
\end{equation}
We prove that 
\begin{equation}
\label{eq:3.25}
w(t)\ge f_k(t),\qquad k=1,2,\dots,
\end{equation}
for almost all $0<\tau<\rho^\theta$ and $\rho^\theta<t<[T-(2\rho)^\theta]/3$. 
By \eqref{eq:3.22} we see that 
\eqref{eq:3.25} holds for $k=1$. 
Assume that \eqref{eq:3.25} holds for some $k\in\{1,2,\dots\}$. 
Then, by \eqref{eq:3.22} we see that 
\begin{equation*}
\begin{split}
w(t) & \ge 2^{-\frac{N}{\theta}}t^{-\frac{N}{\theta}}\int_{\rho^\theta}^t
s^{\frac{N}{\theta}}\left[a_kM_\tau^{p^{k-1}}s^{-\frac{N}{\theta}}\right]^p
\left(\log\frac{s}{\rho^\theta}\right)^{\frac{p(p^{k-1}-1)}{p-1}}\,ds\\
 & =2^{-\frac{N}{\theta}}a_k^pM_\tau^{p^k}t^{-\frac{N}{\theta}}
 \int_{\rho^\theta}^t s^{-1}\left(\log\frac{s}{\rho^\theta}\right)^{\frac{p^k-p}{p-1}}\,ds\\
 & =2^{-\frac{N}{\theta}}a_k^p\,\frac{p-1}{p^k-1}
M_\tau^{p^k}t^{-\frac{N}{\theta}}\left(\log\frac{t}{\rho^\theta}\right)^{\frac{p^k-1}{p-1}}
=f_{k+1}(t)
\end{split}
\end{equation*}
for almost all $0<\tau<\rho^\theta$ and $\rho^\theta<t<[T-(2\rho)^\theta]/3$. 
Therefore we see that \eqref{eq:3.25} holds for all $k=1,2,\dots$. 

On the other hand, 
there exists a constant $\beta>0$ such that 
\begin{equation}
\label{eq:3.26}
a_k\ge\beta^{p^k},\qquad k=1,2,\dots. 
\end{equation}
Indeed, set $b_k:=-p^{-k}\log a_k$. 
It follows from \eqref{eq:3.23} that 
$$
-\log a_{k+1}=-p\log a_k+\log\left[2^{\frac{N}{\theta}}\frac{p^k-1}{p-1}\right],
$$
which yields 
\begin{equation}
\label{eq:3.27}
b_{k+1}-b_k=p^{-k-1}\log\left[2^{\frac{N}{\theta}}\frac{p^k-1}{p-1}\right]
\le p^{-k-1}(Ck+C),
\qquad k=1,2,\dots,
\end{equation}
for some constant $C>0$. 
Since $p>1$, by \eqref{eq:3.27} we see that   
$$
b_{k+1}=b_1+\sum_{i=1}^k (b_{i+1}-b_i)\le C
$$
for $k=1,2,\dots$. 
This implies \eqref{eq:3.26}. 

Taking a sufficiently small $\nu>0$ if necessary, 
by \eqref{eq:3.15}, \eqref{eq:3.24}, \eqref{eq:3.25} and \eqref{eq:3.26} 
we obtain  
\begin{equation*}
\begin{split}
\infty>w(t)\ge f_{k+1}(t)
 & =\biggr[\beta^p M_\tau\log\left(\frac{t}{\rho^\theta}\right)^{\frac{1}{p-1}}\biggr]^{p^k}t^{-\frac{N}{\theta}}
 \log\left(\frac{t}{\rho^\theta}\right)^{-\frac{1}{p-1}}\\
 & \ge\biggr[\beta^p M_\tau\log\left(\frac{T}{5\rho^\theta}\right)^{\frac{1}{p-1}}\biggr]^{p^k} t^{-\frac{N}{\theta}}
 \log\left(\frac{t}{\rho^\theta}\right)^{-\frac{1}{p-1}}, 
\quad k=1,2,\dots,
\end{split}
\end{equation*}
for almost all $T/5<t<T/4$. 
Then it follows that 
$$
\beta^p M_\tau\log\left(\frac{T}{5\rho^\theta}\right)^{\frac{1}{p-1}}\le 1,
$$
which implies the desired inequality 
$$
\int_{B(z,\rho)}u(y,\tau)\,dy
\le\beta^{-p}\biggr[\log\left(\frac{T}{5\rho^\theta}\right)\biggr]^{-\frac{N}{\theta}}
\le C\biggr[\log\left(e+\frac{T}{\rho^\theta}\right)\biggr]^{-\frac{N}{\theta}},
\quad
0<\rho^\theta<\nu T,
$$
for almost all $0<\tau<\rho^\theta$. 
Thus Lemma~\ref{Lemma:3.3} follows.
$\Box$\vspace{5pt}

\noindent
{\bf Proof of Theorem~\ref{Theorem:1.1} in the case $p=p_\theta$.}
By Lemma~\ref{Lemma:2.4} 
we can find a unique Radon measure $\mu$ satisfying \eqref{eq:3.10}. 
Let $\zeta$ be as in \eqref{eq:3.11}. 
Similarly to \eqref{eq:3.12}, 
by Lemma~\ref{Lemma:3.3} we obtain 
\begin{equation*}
\begin{split}
C\biggr[\log\left(e+\frac{T}{\rho^\theta}\right)\biggr]^{-\frac{N}{\theta}}
 & \ge\underset{\tau\to+0}{\mbox{{\rm ess lim}}}\int_{B(z,\rho)}u(y,\tau)\zeta(y-z)\,dy\\
 & =\int_{{\bf R}^N}\zeta(y-z)\,d\mu(y)\ge\mu(B(z,\rho/2))
\end{split}
\end{equation*}
for all $z\in{\bf R}^N$ and $0<\rho<(\nu T)^{1/\theta}$, 
where $\nu$ is as in Lemma~\ref{Lemma:3.3}. 
Then, similarly to \eqref{eq:3.13}, 
we obtain 
$$
\mu(B(z,\sigma))\le 
C\biggr[\log\left(e+\frac{T}{\rho^\theta}\right)\biggr]^{-\frac{N}{\theta}}
$$
for all $z\in{\bf R}^N$ and $0<\sigma<T^{1/\theta}$. 
This is the desired inequality in the case $p=p_\theta$. 
Thus Theorem~\ref{Theorem:1.1} follows. 
$\Box$
\vspace{5pt}

\noindent
{\bf Proof of Corollary~\ref{Corollary:1.1}.}
For almost all $\tau\in(0,T)$, 
$u_\tau(x,t):=u(x,t+\tau)$ is a solution of \eqref{eq:1.5} in ${\bf R}^N\times[0,T-\tau)$ 
with $\mu=u(\tau)$.  
Then, by Lemma~\ref{Lemma:2.4} 
we see that $u(\tau)$ is the initial trace of $u_\tau$. 
Therefore Corollary~\ref{Corollary:1.1} follows from Theorem~\ref{Theorem:1.1} with $\sigma=(T-\tau)^{1/\theta}$.
$\Box$
\vspace{5pt}

\noindent
{\bf Proof of Theorem~\ref{Theorem:1.2}.}
Let $u$ be a solution of \eqref{eq:1.1} in ${\bf R}^N\times(0,T)$, where $0<T<\infty$. 
Let $0<\theta<2$, $0<t<T$ and $n=1,2,\dots$. 
By the Besicovitch covering lemma we can find an integer $m$
depending only on $N$ and a  set $\{x_{k,i}\}_{k=1,\dots,m,\,i\in{\bf N}}\subset{\bf R}^N\setminus B(0,nt^{1/\theta})$ such that 
\begin{equation}
\label{eq:3.28}
B_{k,i}\cap B_{k,j}=\emptyset\quad\mbox{if $i\not=j$}
\qquad\mbox{and}\qquad
{\bf R}^N\setminus B(0,nt^\frac{1}{\theta})\subset\bigcup_{k=1}^m\bigcup_{i=1}^\infty B_{k,i},
\end{equation}
where $B_{k,i}:=B(x_{k,i},t^{1/\theta})$. 
By \eqref{eq:2.1}, \eqref{eq:2.2}, \eqref{eq:3.9} and \eqref{eq:3.28} we obtain  
\begin{equation}
\label{eq:3.29}
\begin{split}
 & \underset{0<\tau<t/2}{\mbox{{\rm ess sup}}}
 \int_{{\bf R}^N\setminus B(0,nt^\frac{1}{\theta})}G(y,t-\tau)u(y,\tau)\,dy\\
 &  \le\sum^m_{k=1}\sum^\infty_{i=1}\,\underset{0<\tau<t/2}{\mbox{{\rm ess sup}}}\int_{B_{k,i}}G(y,t-\tau)u(y,\tau)\,dy\\
 & \le C\,\underset{0<\tau<t/2}{\mbox{{\rm ess sup}}}\sup_{z\in{\bf R}^N}\int_{B(z,t^\frac{1}{\theta})}u(y,\tau)\,dy\\
 & \qquad\qquad
 \times\sum^m_{k=1}\sum^\infty_{i=1}\sup_{0<\tau<t/2}\sup_{y\in B_{k,i}}
 (t-\tau)^{-\frac{N}{\theta}}\left(1+(t-\tau)^{-\frac{1}{\theta}}|y|\right)^{-N-\theta}\\
  & \le Ct^{-\frac{N}{\theta}}\sum^m_{k=1}\sum^\infty_{i=1}\sup_{y\in B_{k,i}}
 \left(1+t^{-\frac{1}{\theta}}|y|\right)^{-N-\theta}.
\end{split}
\end{equation}
On the other hand, 
since
$$
\inf_{y\in B_{k,i}}(1+t^{-\frac{1}{\theta}}|y|)
\ge 1+\frac{1}{4}t^{-\frac{1}{\theta}}\inf_{y\in B_{k,i}}|y|
\ge 1+\frac{1}{4}t^{-\frac{1}{\theta}}(|z|-2t^\frac{1}{\theta})
=\frac{1}{2}+\frac{1}{4}t^{-\frac{1}{\theta}}|z|
$$
for any $z\in B_{k,i}$, 
we have 
$$
\sup_{y\in B_{k,i}}
\left(1+t^{-\frac{1}{\theta}}|y|\right)^{-N-\theta}
\le C\,\,\dashint_{B_{k,i}}
(2+t^{-\frac{1}{\theta}}|z|)^{-N-\theta}\,dz. 
$$
This together with \eqref{eq:3.28} and \eqref{eq:3.29} implies that 
\begin{equation}
\label{eq:3.30}
\begin{split}
 & \underset{0<\tau<t/2}{\mbox{{\rm ess sup}}}\,
 \int_{{\bf R}^N\setminus B(0,nt^\frac{1}{\theta})}G(y,t-\tau)u(y,\tau)\,dy\\
 & \le Ct^{-2\frac{N}{\theta}}\int_{{\bf R}^N\setminus B(0,(n-1)t^{\frac{1}{\theta}})}(2+t^{-\frac{1}{\theta}}|z|)^{-N-\theta}\,dz\\
 & \le Ct^{-\frac{N}{\theta}}\int_{{\bf R}^N\setminus B(0,n-1)}(2+|z|)^{-N-\theta}\,dz\to 0
\end{split}
\end{equation}
as $n\to\infty$. 
Similarly, by Theorem~\ref{Theorem:1.1} we have 
\begin{equation}
\label{eq:3.31}
\begin{split}
 & \int_{{\bf R}^N\setminus B(0,nt^\frac{1}{\theta})}G(y,t)\,d\mu\\
 & \le C\sup_{z\in{\bf R}^N}\mu(B(z,t^\frac{1}{\theta}))
 \sum^m_{k=1}\sum^\infty_{i=1}  \sup_{y\in B_{k,i}}
 t^{-\frac{N}{\theta}}\left(1+t^{-\frac{1}{\theta}}|y|\right)^{-N-\theta}\\
 & \le Ct^{-\frac{N}{\theta}}\int_{{\bf R}^N\setminus B(0,n-1)}(2+|z|)^{-N-\theta}\,dz\to 0
\end{split}
\end{equation}
as $n\to\infty$. 
In particular, 
by Theorem~\ref{Theorem:1.1}, Lemma~\ref{Lemma:2.3}, \eqref{eq:3.30} and \eqref{eq:3.31} 
we see that 
\begin{equation}
\label{eq:3.32}
\int_{{\bf R}^N}G(y,t-\tau)u(y,\tau)\,dy<\infty,
\qquad
\int_{{\bf R}^N} G(y,t)\,d\mu(y)<\infty,
\end{equation}
for almost all $\tau\in(0,t/2)$. 

Let $\eta_n\in C_0({\bf R}^N)$ be such that 
$$
0\le\eta_n\le 1\quad\mbox{in}\quad{\bf R}^N,
\qquad
\eta_n=1\quad\mbox{on}\quad B(0,nt^\frac{1}{\theta}),
\qquad
\eta_n=0\quad\mbox{outside}\quad B(0,2nt^\frac{1}{\theta}). 
$$
It follows from \eqref{eq:3.32} that 
\begin{equation}
\label{eq:3.33}
\begin{split}
 & \left|\int_{{\bf R}^N}G(y,t-\tau)u(y,\tau)\,dy-\int_{{\bf R}^N} G(y,t)\,d\mu(y)\right|\\
 & \le\left|\int_{{\bf R}^N}G(y,t)u(y,\tau)\eta_n(y)\,dy - \int_{{\bf R}^N} G(y,t)\eta_n(y)\,d\mu(y)\right|\\
 & \qquad
 +\left|\int_{{\bf R}^N}[G(y,t-\tau)-G(y,t)]u(y,\tau)\eta_n(y)\,dy\right|\\
 & \qquad\qquad
 +\int_{{\bf R}^N\setminus B(0,nt^\frac{1}{\theta})}G(y,t-\tau)u(y,\tau)\,dy
 +\int_{{\bf R}^N\setminus B(0,nt^\frac{1}{\theta})} G(y,t)\,d\mu(y)
\end{split}
\end{equation}
for $n=1,2,\dots$ and almost all $\tau\in(0,t/2)$. 
By Lemma~\ref{Lemma:2.4} we see that 
\begin{equation}
\label{eq:3.34}
\underset{\tau\to+0}{\mbox{{\rm ess lim}}}\,
\left[\int_{{\bf R}^N}G(y,t)u(y,\tau)\eta_n(y)\,dy - \int_{{\bf R}^N} G(y,t)\eta_n(y)\,d\mu(y)\right]=0.
\end{equation}
Furthermore, by Lemma~\ref{Lemma:2.3} we have 
\begin{equation}
\label{eq:3.35}
\begin{split}
 & \underset{\tau\to+0}{\mbox{{\rm ess limsup}}}\,\left|\int_{{\bf R}^N}[G(y,t-\tau)-G(y,t)]u(y,\tau)\eta_n(y)\,dy\right|\\
 & 
\le \sup_{y\in B(0,2nt^{\frac{1}{\theta}}),s\in(t/2,t)}\,|\partial_t G(y,s)|\,\underset{\tau\to+0}{\mbox{{\rm ess limsup}}}\,
\biggr[\tau\int_{B(0,2nt^\frac{1}{\theta})}u(y,\tau)\,dy\biggr]=0.
\end{split}
\end{equation}
By \eqref{eq:3.33}, \eqref{eq:3.34} and \eqref{eq:3.35} we see that 
\begin{equation*}
\begin{split}
 & \underset{\tau\to+0}{\mbox{{\rm ess limsup}}}\,
 \left|\int_{{\bf R}^N}G(y,t-\tau)u(y,\tau)\,dy - \int_{{\bf R}^N} G(y,t)\,d\mu(y)\right|\\
 & \le
\underset{0<\tau<t/2}{\mbox{{\rm ess sup}}}\int_{{\bf R}^N\setminus B(0,nt^\frac{1}{\theta})}G(y,t-\tau)u(y,\tau)\,dy
 +\int_{{\bf R}^N\setminus B(0,nt^\frac{1}{\theta})} G(y,t)\,d\mu(y)
\end{split}
\end{equation*}
for $n=1,2,\dots$. 
This together with \eqref{eq:3.30} and \eqref{eq:3.31} implies that 
$$
\underset{\tau\to+0}{\mbox{{\rm ess lim}}}\,
\left|\int_{{\bf R}^N}G(y,t-\tau)u(y,\tau)\,dy - \int_{{\bf R}^N} G(y,t)\,d\mu(y)\right|=0.
$$
This together with Definition~\ref{Definition:1.1}~(i) implies that 
$u$ is a solution of \eqref{eq:1.5} in ${\bf R}^N\times[0,T)$. 
Thus Theorem~\ref{Theorem:1.2} follows in the case $0<\theta<2$. 
Similarly, we obtain Theorem~\ref{Theorem:1.2} in the case $\theta=2$, 
and the proof is complete. 
$\Box$
\section{Proof of Theorems~\ref{Theorem:1.3}, \ref{Theorem:1.4} and \ref{Theorem:1.5}.}
We modify the arguments in \cite{IKS} and \cite{RS} 
to prove Theorems~\ref{Theorem:1.3}, \ref{Theorem:1.4} and \ref{Theorem:1.5}. 
In the rest of this paper, 
for any two nonnegative functions
$f_1$ and $f_2$ defined in a subset $D$ of $[0,\infty)$,
we write
$f_1(t)\asymp f_2(t)$ for all $t\in D$ 
if $C^{-1}f_2(t)\le f_1(t)\le Cf_2(t)$ for all $t\in D$. 
\vspace{5pt}

\noindent
{\bf Proof of Theorem~\ref{Theorem:1.3}.}
It suffices to consider the case $T=1$. 
Indeed, for any solution $u$ of \eqref{eq:1.5} in ${\bf R}^N\times[0,T)$, where $0<T<\infty$, 
we see that
$u_\lambda(x,t):=\lambda^{\theta/(p-1)}u(\lambda x,\lambda^\theta t)$
with $\lambda:=T^{1/\theta}$ is a solution of \eqref{eq:1.5} in ${\bf R}^N\times[0,1)$. 

Set $w(x,t):=S(t)\mu$. Then it follows from Lemma~\ref{Lemma:2.1} that 
\begin{equation*}
\begin{split}
F[w](t):= & S(t)\mu+\int_0^tS(t-s)(2w(s))^p\,ds
\le w(t)+2^p w(t)\int_0^t\|w(s)\|_{L^\infty({\bf R}^N)}^{p-1}\,ds\\
 \le & w(t)+Cw(t)\int_0^t s^{-\frac{N(p-1)}{\theta}}\left[\sup_{x\in{\bf R}^N}\mu(B(x,s^{\frac{1}{\theta}}))\right]^{p-1}\,ds. 
\end{split}
\end{equation*}
This together with \eqref{eq:1.9} and $1<p<p_\theta$ implies that 
$$
F[w](t)\le w(t)+C\gamma_2^{p-1}t^{1-\frac{N(p-1)}{\theta}}w(t)
\le [1+C\gamma_2^{p-1}]w(t)
$$
for all $0\le t<1$. 
Therefore, taking a sufficiently small $\gamma_2>0$ if necessary, 
we obtain 
$F[w](t)\le 2w(t)$ for $0\le t<1$. 
This means that $2w(t)$ is a supersolution of \eqref{eq:1.1} in ${\bf R}^N\times[0,1)$. 
Then Theorem~\ref{Theorem:1.3} follows from Lemma~\ref{Lemma:2.2}.
$\Box$\vspace{5pt}

\noindent
{\bf Proof of Theorem~\ref{Theorem:1.4}.}
Similarly to the proof of Theorem~\ref{Theorem:1.3}, 
it suffices to consider the case $T=1$. 
Let $1<\alpha<p$ and set $w(x,t):=[S(t)\mu^\alpha]^{1/\alpha}$. 
It follows from the Jensen inequality that 
\begin{equation*}
\begin{split}
F[w](t):= & S(t)\mu+\int_0^tS(t-s)(2w(s))^p\,ds
\le w(t)+2^p[S(t)\mu^\alpha]\int_0^t\|S(t)\mu^\alpha\|_{L^\infty({\bf R}^N)}^{\frac{p}{\alpha}-1}\,ds\\
 \le & w(t)+Cw(t)\|S(t)\mu^\alpha\|_{L^\infty({\bf R}^N)}^{1-\frac{1}{\alpha}}
 \int_0^t \|S(s)\mu^\alpha\|_{L^\infty({\bf R}^N)}^{\frac{p}{\alpha}-1}\,ds
\end{split}
\end{equation*}
for all $0\le t<1$. 
Furthermore, by Lemma~\ref{Lemma:2.1} and \eqref{eq:1.10} we have
$$
\|S(t)\mu^\alpha\|_{L^\infty({\bf R}^N)}
\le C\gamma_3^\alpha t^{-\frac{\alpha}{p-1}},\qquad 0\le t<1. 
$$
This implies that 
$$
F[w](t)\le w(t)+C\gamma_3^{p-1}w(t)t^{-\frac{\alpha-1}{p-1}}\int_0^t s^{-\frac{p-\alpha}{p-1}}\,ds
\le [1+C\gamma_3^{p-1}]w(t),
\quad 0\le t<1.
$$
Therefore, taking a sufficiently small $\gamma_3>0$ if necessary, 
we obtain 
$F[w](t)\le 2w(t)$ for $0\le t<1$. 
Then, similarly to Theorem~\ref{Theorem:1.3}, 
Theorem~\ref{Theorem:1.4} follows from Lemma~\ref{Lemma:2.2}.
$\Box$\vspace{5pt}

\noindent
{\bf Proof of Theorem~\ref{Theorem:1.5}.}
Similarly to the proof of Theorem~\ref{Theorem:1.3}, 
it suffices to consider the case $T=1$. 
Let $\beta>0$ and $\rho=\rho(s)$ be as in \eqref{eq:1.11}. 
Let $L\ge e$ be such that  
\begin{itemize}
  \item[{\rm (a)}] 
  $\Psi_{\beta,L}(s):=s[\log (L+s)]^\beta$ is positive and convex in $(0,\infty)$;
  \item[{\rm (b)}] 
  $s^p/\Psi_{\beta,L}(s)$ and $\Psi_{\beta,L}(s)/s$ are monotone increasing in $(0,\infty)$.
\end{itemize}
Let $\mu$ be a nonnegative measurable function in ${\bf R}^N$ satisfying \eqref{eq:1.12}. 
Since $\Psi_\beta(\tau)\asymp \Psi_{\beta,L}(\tau)$ for $\tau>0$,  
it follows that 
\begin{equation}
\label{eq:4.1}
\Psi_{\beta,L}^{-1}\left[\,\dashint_{B(x,\sigma)}\Psi_{\beta,L}(\mu(y))\,dy\,\right]\le C\gamma\rho(\sigma)
\end{equation}
for all $x\in{\bf R}^N$ and $0<\sigma<1$.
Set 
$$
w(x,t):=\Psi_{\beta,L}^{-1}\left[\,S(t)\Psi_{\beta,L}(\mu)\,\right]. 
$$
By \eqref{eq:4.1}
we apply Lemma~\ref{Lemma:2.1} to obtain 
$$
\|S(t)\Psi_{\beta,L}(\mu)\|_{L^\infty({\bf R}^N)}
=\|\Psi_{\beta,L}(w(t))\|_{L^\infty({\bf R}^N)}\le C\Psi_{\beta,L}(C\gamma\rho(t^\frac{1}{\theta}))
\le C\Psi_{\beta,L}(\gamma\rho(t^\frac{1}{\theta})),
$$
which implies that 
\begin{equation}
\label{eq:4.2}
\|w(t)\|_{L^\infty({\bf R}^N)}\le\Psi_{\beta,L}^{-1}\biggr[C\Psi_{\beta,L}(\gamma\rho(t^\frac{1}{\theta}))\biggr]
\end{equation}
for $0<t<1$. 
Define
$$
F[w](t):=S(t)\mu+\int_0^tS(t-s)(2w(s))^p\,ds,\qquad t>0.
$$
It follows from the Jensen inequality and property~(a) that 
$S(t)\mu\le w(t)$. Then 
\begin{equation}
\label{eq:4.3}
\begin{split}
F[w](t) 
 & \le w(t)+2^p\int_0^t S(t-s)\left[\frac{w(s)^p}{S(s)\Psi_{\beta,L}(\mu)}S(s)\Psi_{\beta,L}(\mu)\right]ds\\
 & \le w(t)+2^p\biggr[\,\int_0^t \left\|\frac{w(s)^p}{S(s)\Psi_{\beta,L}(\mu)}\right\|_{L^\infty({\bf R}^N)}\,ds\,\biggr]S(t)\Psi_{\beta,L}(\mu)\\
 & \le w(t)+C\biggr[\,\int_0^t \left\|\frac{w(s)^p}{\Psi_{\beta,L}(w(s))}\right\|_{L^\infty({\bf R}^N)}\,ds\,\biggr]
 \left\|\frac{\Psi_{\beta,L}(w(t))}{w(t)}\right\|_{L^\infty({\bf R}^N)} w(t)
\end{split}
\end{equation} 
for $t>0$. 
On the other hand, 
by property~(b) and \eqref{eq:4.2} we see that  
\begin{equation}
\label{eq:4.4}
\left\|\frac{w(s)^p}{\Psi_{\beta,L}(w(s))}\right\|_{L^\infty({\bf R}^N)}
\le\frac{\|w(s)\|_{L^\infty({\bf R}^N)}^p}{\Psi_{\beta,L}(\|w(s)\|_{L^\infty({\bf R}^N)})}
\le\frac{[{\Psi_{\beta,L}^{-1}}(C\Psi_{\beta,L}(\gamma\rho(s^\frac{1}{\theta})))]^p}{C\Psi_{\beta,L}(\gamma\rho(s^\frac{1}{\theta}))}
\end{equation}
for $0<s<1$.
By \eqref{eq:1.11} we have
$$
\Psi_{\beta,L}(\gamma\rho(s^\frac{1}{\theta}))=\gamma\rho(s^\frac{1}{\theta})[\log(L+\gamma\rho(s^\frac{1}{\theta}))]^\beta
\asymp
\gamma s^{-\frac{N}{\theta}}\biggr[\log\biggr(e+\frac{1}{s}\biggr)\biggr]^{-\frac{N}{\theta}+\beta}
$$
for all $s\in(0,1)$. 
Since 
$\Psi_{\beta,L}^{-1}(\tau)\asymp \tau[\log(e+\tau)]^{-\beta}$ for all $\tau>0$, 
it follows that 
$$
\Psi_{\beta,L}^{-1}(C\Psi_{\beta,L}(\gamma\rho(s^\frac{1}{\theta})))\asymp \gamma s^{-\frac{N}{\theta}}
\biggr[\log\biggr(e+\frac{1}{s}\biggr)\biggr]^{-\frac{N}{\theta}}
$$
for all $s\in(0,1)$. These together with \eqref{eq:4.4} imply that 
\begin{equation}
\label{eq:4.5}
\begin{split}
\left\|\frac{w(s)^p}{\Psi_{\beta,L}(w(s))}\right\|_{L^\infty({\bf R}^N)}
\le C\gamma^{\frac{\theta}{N}}s^{-1}\biggr[\log\biggr(e+\frac{1}{s}\biggr)\biggr]^{-1-\beta}
\end{split}
\end{equation}
for all $s\in(0,1)$. Similarly, by \eqref{eq:4.2} and property~(b) we have
\begin{equation}
\label{eq:4.6}
\left\|\frac{\Psi_{\beta,L}(w(t))}{w(t)}\right\|_{L^\infty({\bf R}^N)}
\le\frac{C\Psi_{\beta,L}(\gamma\rho(t^\frac{1}{\theta}))}{\Psi_{\beta,L}^{-1}(C\Psi_{\beta,L}(\gamma\rho(t^\frac{1}{\theta})))}
\le C\biggr[\log\biggr(e+\frac{1}{t}\biggr)\biggr]^\beta
\end{equation}
for all $t\in(0,1)$.
By \eqref{eq:4.5} and \eqref{eq:4.6} we obtain 
\begin{equation}
\label{eq:4.7}
\begin{split}
 & \biggr[\,\int_0^t \left\|\frac{w(s)^p}{\Psi_{\beta,L}(w(s))}\right\|_{L^\infty({\bf R}^N)}\,ds\,\biggr]
\left\|\frac{\Psi_{\beta,L}(w)}{w(t)}\right\|_{L^\infty({\bf R}^N)}\\
 & \le C\gamma^{\frac{\theta}{N}}\biggr[\log\biggr(e+\frac{1}{t}\biggr)\biggr]^\beta
\int_0^t s^{-1}\biggr[\log\biggr(e+\frac{1}{s}\biggr)\biggr]^{-1-\beta}\,ds
\le C\gamma^{\frac{\theta}{N}}
\end{split}
\end{equation}
for all $0\le t<1$. 
Therefore, taking a sufficiently small $\gamma$ if necessary, 
we deduce from \eqref{eq:4.3} and \eqref{eq:4.7} that 
$F[w](t)\le 2w(t)$ for $0\le t<1$. 
Then, similarly to Theorem~\ref{Theorem:1.3}, 
Theorem~\ref{Theorem:1.5} follows from Lemma~\ref{Lemma:2.2}. 
$\Box$\vspace{5pt}

As a corollary of Theorem~\ref{Theorem:1.5}, we have
\begin{corollary}
\label{Corollary:4.1}
Assume the same conditions as in Theorem~{\rm\ref{Theorem:1.5}}. 
Let $\alpha>1$ and $\rho=\rho(s)$ be as in \eqref{eq:1.11}. 
Then there exists $\gamma>0$ such that, 
if $\mu$ is a nonnegative measurable function in ${\bf R}^N$ satisfying 
$$
\sup_{x\in{\bf R}^N}\left[\,\dashint_{B(x,\sigma)}
(T^\frac{1}{p-1}\mu(y))^\alpha\,dy\,\right]^{\frac{1}{\alpha}}\le\gamma\rho(\sigma T^{-\frac{1}{\theta}}),
\qquad
0<\sigma<T^{\frac{1}{\theta}},
$$
for some $T>0$, 
then problem~\eqref{eq:1.5} possesses a solution in ${\bf R}^N\times[0,T)$. 
\end{corollary}
{\bf Proof.}
Let $T>0$, $\alpha>1$ and $\beta>0$. 
Let $\Phi=\Phi(s)$ be a nonnegative and convex function in $[0,\infty)$ such that 
$C^{-1}s^\alpha\le \Phi(\Psi_\beta(s))\le Cs^\alpha$ in $[0,\infty)$ for some positive constant $C$. 
Then it follows from the Jensen inequality that 
\begin{equation*}
\begin{split}
\Psi_\beta^{-1}\left[\,\dashint_{B(x,\sigma)}
\Psi_\beta(T^\frac{1}{p-1}\mu(y))\,dy\right]
 & \le\Psi_\beta^{-1}\left[\Phi^{-1}
\left[\,\dashint_{B(x,\sigma)}
\Phi\left(\Psi_\beta(T^\frac{1}{p-1}\mu(y))\right)\,dy\right]\right]\\
 & \le C\left[\,\dashint_{B(x,\sigma)}
(T^\frac{1}{p-1}\mu(y))^\alpha\,dy\right]^{\frac{1}{\alpha}}
\end{split}
\end{equation*}
for $0<\sigma<T^{1/\theta}$. 
Therefore Corollary~\ref{Corollary:4.1} follows from Theorem~{\rm\ref{Theorem:1.5}}. 
$\Box$\vspace{5pt}
\newline
Corollary~\ref{Corollary:4.1} is also a refinement of Theorem~\ref{Theorem:1.4}. 
\section{Application}
Since the minimal solution is unique, 
we can define the maximal existence time $T(\mu)$  of the minimal solution of \eqref{eq:1.5}. 
See also Lemma~\ref{Lemma:2.2}. 

For problem~\eqref{eq:1.5} with $\theta=2$, 
Lee and Ni~\cite{LN} obtained optimal estimates of $T(\lambda\phi)$ as $\lambda\to+0$ 
by use of the behavior of $\phi$ at the space infinity. 
In this section, as an application of our theorems, 
we show that similar estimates of $T(\lambda\phi)$ as in \cite{LN} 
hold in the case $0<\theta<2$. 
Theorems~\ref{Theorem:5.1} and \ref{Theorem:5.2} 
are generalizations of \cite[Theorem~3.15 and 3.21]{LN}, respectively.
\begin{theorem}
\label{Theorem:5.1}
Let $N\ge 1$, $0<\theta\le 2$ and $p>1$. 
Let $A>0$ and $\phi$ be a nonnegative measurable function in ${\bf R}^N$ such that 
\begin{equation}
\label{eq:5.1}
0\le\phi(x)\le(1+|x|)^{-A}
\end{equation}
for almost all $x\in{\bf R}^N$. 
\begin{itemize}
  \item[{\rm (i)}] 
  Let $p=p_\theta$ and $A\ge \theta/(p-1)=N$. 
  Then there exists a positive constant $C_1$ such that 
  $$
  \log T(\lambda\phi)\ge
  \left\{
  \begin{array}{ll}
  C_1\lambda^{-(p-1)} & \mbox{if}\quad A>N,\\
  C_1\lambda^{-\frac{p-1}{p}} & \mbox{if}\quad A=N,\\
  \end{array}
  \right.
  $$
  for all small enough $\lambda>0$. 
  \item[{\rm (ii)}] 
  Let $1<p<p_\theta$ or  $A<\theta/(p-1)$. 
  Then there exists a positive constant $C_2$ such that 
  $$
  T(\lambda\phi)\ge
  \left\{
  \begin{array}{ll}
  C_2\lambda^{-\left(\frac{1}{p-1}-\frac{1}{\theta} \min \{A,N\} \right)^{-1}} & \mbox{if}\quad A\not=N,\vspace{3pt}\\
  C_2\left(\frac{\lambda^{-1}}{\log(\lambda^{-1})}\right)^{\left(\frac{1}{p-1}-\frac{N}{\theta}\right)^{-1}} & \mbox{if}\quad A=N,\\
  \end{array}
  \right.
  $$
  for all small enough $\lambda>0$. 
\end{itemize}
\end{theorem}
{\bf Proof.}
We apply Theorem~\ref{Theorem:1.5} with 
and prove assertion~(i). 
Let $p=p_\theta$ and set 
\begin{equation}
\label{eq:5.2}
\Psi(s):=s[\log(e+s)]^{\frac{N}{\theta}},
\qquad
\rho(s):=
s^{-N}\biggr[\log\biggr(e+\frac{1}{s}\biggr)\biggr]^{-\frac{N}{\theta}},
\end{equation}
for $s>0$ (see \eqref{eq:1.11}).  
For any $0<\epsilon<1$ and $0<\delta<1$, set 
\begin{equation}
\label{eq:5.3}
T_\lambda:=
\left\{
\begin{array}{ll}
\exp(\delta \lambda^{-(p-1)}) & \mbox{if}\quad A>N,\\
\exp(\delta \lambda^{-\frac{p-1}{p}}) & \mbox{if}\quad A=N.\\
\end{array}
\right.
\end{equation}
By \eqref{eq:5.1}, \eqref{eq:5.2} and \eqref{eq:5.3}
we have 
\begin{equation}
\label{eq:5.4}
\begin{split}
 & \dashint_{B(x,\sigma)}\Psi(T_\lambda^\frac{1}{p-1} \lambda \phi(y) )\,dy
 =\dashint_{B(x,\sigma)}T_\lambda^\frac{1}{p-1} \lambda \phi(y) 
 [\log(e+T_\lambda^\frac{1}{p-1} \lambda \phi(y))]^{\frac{N}{\theta}}\,dy\\
 & \le C\sigma^{-N} \lambda T_\lambda^\frac{1}{p-1}[\log(e+\lambda T_\lambda^\frac{1}{p-1})]^{\frac{N}{\theta}}
 \int_{B(0,\sigma)}(1+|y|)^{-A}\,dy\\
 & \le\left\{
 \begin{array}{ll}
 C \lambda T_\lambda^\frac{1}{p-1}[\log(e+\lambda T_\lambda^\frac{1}{p-1})]^{\frac{N}{\theta}} & \mbox{if}\quad 0<\sigma\le 1,\\
 C\delta^{\frac{N}{\theta}}\sigma^{-N}
 T_\lambda^\frac{1}{p-1} & \mbox{if}\quad 1<\sigma<T_\lambda^{\frac{1}{\theta}},\,\,A>N,\\
 C\delta^{\frac{N}{\theta}}\sigma^{-N}\log(e+\sigma)
 \lambda^{1-\frac{1}{p}} T_\lambda^\frac{1}{p-1} & \mbox{if}\quad 1<\sigma<T_\lambda^{\frac{1}{\theta}},\,\,A=N,\\
 \end{array}
 \right.\\
 & \le\left\{
 \begin{array}{ll}
 C \lambda T_\lambda^\frac{1}{p-1}[\log(e+\lambda T_\lambda^\frac{1}{p-1})]^{\frac{N}{\theta}} & \mbox{if}\quad 0<\sigma\le 1,\\
 C\delta^{\frac{N}{\theta}}\sigma^{-N}
 T_\lambda^\frac{1}{p-1} & \mbox{if}\quad 1<\sigma<T_\lambda^{\frac{1}{\theta}},
 \end{array}
 \right.
\end{split}
\end{equation}
for all $x\in{\bf R}^N$, $0<\sigma<T_\lambda^{1/\theta}$ and all small enough $\lambda>0$. 
Then we have 
\begin{equation}
\label{eq:5.5}
\sup_{x\in{\bf R}^N}\Psi^{-1}\biggr(\dashint_{B(x,\sigma)}\Psi(T_\lambda^\frac{1}{p-1} \lambda \phi(y) )\,dy\biggr)
\le C \lambda T_\lambda^{\frac{1}{p-1}},
\qquad 0<\sigma<1, 
\end{equation}
for all small enough $\lambda>0$. 
On the other hand, let $L\ge e$ be such that 
\begin{equation}
\label{eq:5.6}
\mbox{$s\log[L+s]^{-N/\theta}$ is monotone in $[0,\infty)$}. 
\end{equation}
Since $\Psi(s)^{-1}\asymp s[\log(e+s)]^{-N/\theta}\asymp s\log[L+s]^{-N/\theta}$ for all $s>0$, 
it follows from \eqref{eq:5.4} and \eqref{eq:5.6} that 
\begin{equation}
\label{eq:5.7}
\begin{split}
 & \sup_{x\in{\bf R}^N}\Psi^{-1}\biggr(\dashint_{B(x,\sigma)}\Psi(T_\lambda^\frac{1}{p-1} \lambda \phi(y) )\,dy\biggr)
\le C\delta^{\frac{N}{\theta}}\sigma^{-N} 
 T_\lambda^\frac{1}{p-1}
 \displaystyle{\biggr[\log\biggr(L+\delta^{\frac{N}{\theta}}\frac{T^{\frac{N}{\theta}}_\lambda}{\sigma^N}\biggr)\biggr]^{-\frac{N}{\theta}}}\\
 & \qquad\quad
 \le C\delta^{\frac{N}{\theta}}\sigma^{-N} 
 T_\lambda^\frac{1}{p-1}
 \displaystyle{\biggr[\log\biggr(L+\delta^{\frac{1}{2\theta}}\frac{T^{\frac{1}{2\theta}}_\lambda}{\sigma^{\frac{1}{2}}}\biggr)\biggr]^{-\frac{N}{\theta}}}\\
 & \qquad\quad
 =C\delta^{\frac{N}{\theta}}\sigma^{-N} 
 T_\lambda^\frac{1}{p-1}
 \left(\delta^{\frac{1}{2\theta}}\frac{T^{\frac{1}{2\theta}}_\lambda}{\sigma^{\frac{1}{2}}}\right)^{-1}
 \delta^{\frac{1}{2\theta}}\frac{T^{\frac{1}{2\theta}}_\lambda}{\sigma^{\frac{1}{2}}}
 \displaystyle{\biggr[\log\biggr(L+\delta^{\frac{1}{2\theta}}\frac{T^{\frac{1}{2\theta}}_\lambda}{\sigma^{\frac{1}{2}}}\biggr)\biggr]^{-\frac{N}{\theta}}}\\
 & \qquad\quad
 \le C\delta^{\frac{N}{\theta}}\sigma^{-N} 
 T_\lambda^\frac{1}{p-1}
 \left(\delta^{\frac{1}{2\theta}}\frac{T^{\frac{1}{2\theta}}_\lambda}{\sigma^{\frac{1}{2}}}\right)^{-1}
 \frac{T^{\frac{1}{2\theta}}_\lambda}{\sigma^{\frac{1}{2}}}
 \displaystyle{\biggr[\log\biggr(L+\frac{T^{\frac{1}{2\theta}}_\lambda}{\sigma^{\frac{1}{2}}}\biggr)\biggr]^{-\frac{N}{\theta}}}\\
 & \qquad\quad
 \le C\delta^{\frac{N}{\theta}-\frac{1}{2\theta}}\sigma^{-N} T_\lambda^\frac{1}{p-1}
 \displaystyle{\biggr[\log\biggr(L+\frac{T^{\frac{1}{\theta}}_\lambda}{\sigma}\biggr)\biggr]^{-\frac{N}{\theta}}}
\end{split}
\end{equation}
for all $\sigma\ge 1$ and small enough $\lambda>0$.
Taking a sufficiently small $\delta>0$ if necessary, 
by \eqref{eq:5.2}, \eqref{eq:5.5} and \eqref{eq:5.7} 
we obtain 
$$
\sup_{x\in{\bf R}^N}\Psi^{-1}\biggr(\dashint_{B(x,\sigma)}\Psi(T_\lambda^\frac{1}{p-1} \lambda \phi(y) )\,dy\biggr)
\le\gamma_4\rho(\sigma T_\lambda^{-\frac{1}{\theta}}),
\quad
0<\sigma<T_\lambda^{\frac{1}{\theta}},
$$
for all small enough $\lambda>0$, 
where $\gamma_4$ is the constant given in Theorem~\ref{Theorem:1.5} with $\beta=N/\theta$. 
Then, by Theorem~\ref{Theorem:1.5} 
we see that problem~\eqref{eq:1.5} with $\mu=\lambda\phi$ possesses a solution in ${\bf R}^N\times[0,T_\lambda)$. 
This implies that $T(\lambda\phi)\ge T_\lambda$ for all small enough $\lambda>0$ 
and assertion~(i) follows. 

We prove assertion~(ii) in the case where $p\ge p_\theta$ and $A<\theta/(p-1)$. 
It follows that $A<\theta/(p-1)\le N$. 
For $0<\lambda<1$ and $0<\delta<1$, set 
$$
\tilde{T}_\lambda:=\delta\lambda^{-\left(\frac{1}{p-1}-\frac{A}{\theta}\right)^{-1}}. 
$$
Let $1<\alpha<p$ be such that $A\alpha<N$. 
Then  
\begin{equation}
\label{eq:5.8}
\begin{split}
\left(\dashint_{B(x,\sigma)}(\lambda \phi(y))^\alpha\,dy\right)^{\frac{1}{\alpha}}
 & \le \lambda
 \left(\dashint_{B(x,\sigma)}(1+|y|)^{-A\alpha}\,dy\right)^{\frac{1}{\alpha}}
 \le C\lambda\sigma^{-A}
\end{split}
\end{equation}
for all $x\in{\bf R}^N$ and $0<\sigma<\tilde{T}_\lambda^{1/\theta}$. 
On the other hand, 
it follows that 
\begin{equation}
\label{eq:5.9}
\lambda\sigma^{-A}
=\sigma^{-\frac{\theta}{p-1}}\cdot \lambda\sigma^{\frac{\theta}{p-1}-A}
\le\sigma^{-\frac{\theta}{p-1}}\cdot \lambda \tilde{T}_\lambda^{\frac{1}{p-1}-\frac{A}{\theta}}
=\delta^{\frac{1}{p-1}-\frac{A}{\theta}}\sigma^{-\frac{\theta}{p-1}}
\end{equation}
for $0<\sigma<\tilde{T}_\lambda^{1/\theta}$. 
Then, taking a sufficiently small $\delta>0$ if necessary, 
by \eqref{eq:5.8} and \eqref{eq:5.9} 
we obtain \eqref{eq:1.10} for all $0<\sigma<\tilde{T}_\lambda^{1/\theta}$.
Therefore, similarly to assertion~(i), 
we deduce from Theorem~\ref{Theorem:1.4} that $T(\lambda\phi)\ge \tilde{T}_\lambda$ for $\lambda\in(0,1)$.  
Thus assertion~(ii) follows in the case where $p\ge p_\theta$ and $A<\theta/(p-1)$.

It remains to prove assertion~(ii) in the case $1<p<p_\theta$. 
For $0<\lambda<1$ and $0<\delta<1$, 
set 
$$
\hat{T}_\lambda:=
\left\{
\begin{array}{ll}
\delta\lambda^{-\left(\frac{1}{p-1}-\frac{1}{\theta} \min \{A,N\}\right)^{-1}} & \mbox{if}\quad A\not=N,\vspace{3pt}\\
\delta\left(\frac{\lambda^{-1}}{\log(\lambda^{-1})}\right)^{\left(\frac{1}{p-1}-\frac{N}{\theta}\right)^{-1}} & \mbox{if}\quad A=N.\\
\end{array}
\right.
$$ 
It follows that
$$
\sup_{x\in{\bf R}^N}\int_{B(x,\hat{T}_\lambda^\frac{1}{\theta})}\lambda \phi(y)\, dy \le
\left\{
\begin{array}{ll}
C\lambda & \mbox{if}\quad A>N,\vspace{3pt}\\
C\lambda \log(1+\hat{T}_\lambda^\frac{1}{\theta}) & \mbox{if}\quad A=N,\vspace{3pt}\\
C\lambda \hat{T}_\lambda^\frac{N-A}{\theta}& \mbox{if}\quad 0<A<N.
\end{array}
\right.
$$
Then, taking a sufficiently small $\delta>0$ if necessary, 
we obtain 
$$
\sup_{x\in{\bf R}^N}\int_{B(x,\hat{T}_\lambda^\frac{1}{\theta})}\lambda \phi(y)\, dy 
\le \gamma_2 \hat{T}_\lambda^{\frac{N}{\theta}-\frac{1}{p-1}}
$$
for all small enough $\lambda>0$, 
where $\gamma_2$ is as in Theorem~\ref{Theorem:1.3}. 
Therefore, by Theorem~\ref{Theorem:1.3} 
we see that $T(\lambda\phi)\ge \hat{T}_\lambda$ for all small enough $\lambda>0$, 
and assertion~(ii) follows in the case $1<p<p_\theta$. 
Thus the proof of Theorem~\ref{Theorem:5.1} is complete. 
$\Box$
\begin{theorem}
\label{Theorem:5.2}
Let $N\ge 1$, $0<\theta\le 2$ and $p>1$. 
Let $A>0$ and $\phi$ be a nonnegative $L^\infty({\bf R}^N)$-function such that 
\begin{equation}
\label{eq:5.10}
\phi(x)\ge(1+|x|)^{-A}
\end{equation}
for almost all $x\in{\bf R}^N$. 
\begin{itemize}
  \item[{\rm (i)}] 
  Let $p=p_\theta$ and $A\ge \theta/(p-1)=N$. 
  Then there exists a positive constant $C_1$ such that 
  \begin{equation}
  \label{eq:5.11}
  \log T(\lambda\phi)\le
  \left\{
  \begin{array}{ll}
  C_1\lambda^{-(p-1)} & \mbox{if}\quad A>N,\\
  C_1\lambda^{-\frac{p-1}{p}} & \mbox{if}\quad A=N,\\
  \end{array}
  \right.
  \end{equation}
  for all small enough $\lambda>0$. 
  \item[{\rm (ii)}] 
  Let $1<p<p_\theta$ or  $A<\theta/(p-1)$. 
  Then there exists a positive constant $C_2$ such that 
  \begin{equation}
  \label{eq:5.12}
  T(\lambda\phi)\le
  \left\{
  \begin{array}{ll}
  C_2\lambda^{-\left(\frac{1}{p-1}-\frac{1}{\theta} \min \{A,N\} \right)^{-1}} & \mbox{if}\quad A\not=N,\vspace{3pt}\\
  C_2\left(\frac{\lambda^{-1}}{\log(\lambda^{-1})}\right)^{\left(\frac{1}{p-1}-\frac{N}{\theta}\right)^{-1}} & \mbox{if}\quad A=N,\\
  \end{array}
  \right.
  \end{equation}
  for all small enough $\lambda>0$. 
\end{itemize}
\end{theorem}
{\bf Proof.}
Since $\phi\in L^\infty({\bf R}^N)$, 
by Theorem~\ref{Theorem:1.4} 
we have 
\begin{equation}
\label{eq:5.13}
T(\lambda\phi)\ge C\lambda^{-(p-1)}
\end{equation}
for all small enough $\lambda>0$. 
It follows from \eqref{eq:5.10} that 
\begin{equation}
\label{eq:5.14}
\int_{B(0,\sigma)}\lambda \phi(y)\,dy \ge
\lambda\int_{B(0,\sigma)}(1+|y|)^{-A}\,dy\ge
\left\{
\begin{array}{ll}
C\lambda & \mbox{if}\quad \sigma>1,\,\,A>N,\vspace{3pt}\\
C\lambda\log(e+\sigma) & \mbox{if}\quad \sigma>1,\,\,A=N,\vspace{3pt}\\
C\lambda\sigma^{N-A} & \mbox{if}\quad \sigma>1,\,\,A<N. 
\end{array}
\right.
\end{equation}
On the other hand, in the case $p=p_\theta$, 
by Theorem~\ref{Theorem:1.1}
we can find a constant $\gamma_1$ such that 
$$
\int_{B(0,\sigma)}\lambda \phi(y)\,dy\le \gamma_1\,
\biggr[\log\biggl(e+\frac{T(\lambda\phi)^\frac{1}{\theta}}{\sigma}\biggr)\biggr]^{-\frac{N}{\theta}}
$$
for all $0<\sigma<T(\lambda\phi)^{1/\theta}$ and $\lambda>0$. 
This together with \eqref{eq:5.13} implies that 
\begin{eqnarray}
\label{eq:5.15}
 & & \int_{B(0,T(\lambda\phi)^{\frac{1}{2\theta}})}\lambda \phi(y)\,dy\le 
C\gamma_1[\log T(\lambda\phi)]^{-\frac{N}{\theta}},\\
\label{eq:5.16}
 & & \int_{B(0,T(\lambda\phi)^{\frac{1}{\theta}}/2)}\lambda \phi(y)\,dy\le 
C\gamma_1,
\end{eqnarray}
for all small enough $\lambda>0$. 
Setting $\sigma=T(\lambda\phi)^{1/2\theta}\ge 1$, 
by \eqref{eq:5.14} and \eqref{eq:5.15} 
we obtain \eqref{eq:5.11}. 
Thus assertion~(i) follows. 
Furthermore, 
setting $\sigma=T(\lambda\phi)^{1/\theta}/2\ge 1$, 
by \eqref{eq:5.14} and \eqref{eq:5.16}  
we obtain \eqref{eq:5.12} in the case where $p=p_\theta$ and $A<N$. 

We prove assertion~(ii) in the case $p<p_\theta$.
By Theorem~\ref{Theorem:1.1} with $1<p<p_\theta$ 
we have 
\begin{equation}
\label{eq:5.17}
\int_{B(0,\sigma)}\lambda \phi(y)\,dy\le \gamma_1\,
T(\lambda\phi)^{\frac{N}{\theta}-\frac{1}{p-1}}
\end{equation}
for all $0<\sigma<T(\lambda\phi)^{1/\theta}$ and $\lambda>0$. 
Then, by \eqref{eq:5.13}, \eqref{eq:5.14} and \eqref{eq:5.17} with $\sigma=T(\lambda\phi)^{1/\theta}/2\ge 1$ 
we obtain \eqref{eq:5.12} in the case $1<p<p_\theta$. 
Similarly, we obtain \eqref{eq:5.12} in the case $p>p_\theta$ and $A<\theta/(p-1)$. 
Thus assertion~(ii) follows, 
and the proof of Theorem~\ref{Theorem:5.2} is complete. 
$\Box$
\vspace{7pt}

\noindent
{\bf Acknowledgment.} 
The second author of this paper was supported
by the Grant-in-Aid for Scientific Research (A)(No.~15H02058),
from Japan Society for the Promotion of Science.
\bibliographystyle{amsplain}

\bigskip
\noindent Addresses:

\smallskip
\noindent K. H.:  Mathematical Institute, Tohoku University,
Aoba, Sendai 980-8578, Japan\\
\noindent 
E-mail: {\tt kotaro.hisa.s5@dc.tohoku.ac.jp}\\

\noindent K. I.: Mathematical Institute, Tohoku University,
Aoba, Sendai 980-8578, Japan\\
\noindent 
E-mail: {\tt ishige@math.tohoku.ac.jp}\\
\end{document}